\documentclass[10pt]{article}  
\usepackage{amssymb}  
\usepackage{amsthm}
\usepackage{amsmath}
\usepackage{graphicx}
 \usepackage{mathtools}
 
 \usepackage{mathrsfs}

\usepackage[alphabetic,lite]{amsrefs}

\newcommand{\CC}{\mathbb C}
\newcommand{\FF}{\mathbb F}
\newcommand{\bF}{\mathbb F}

\newcommand{\ZZ}{\mathbb Z}
\newcommand{\bZ}{\mathbb Z}

\newcommand{\QQ}{\mathbb Q}
\newcommand{\RR}{\mathbb R}

\newcommand{\frakf}{\mathfrak f}
\newcommand{\frakg}{\mathfrak g}

\newcommand{\la}{\langle}
\newcommand{\ra}{\rangle}

\newcommand{\al}{\alpha}
\newcommand{\be}{\beta}
\newcommand{\del}{\delta}
\newcommand{\lam}{\lambda}
\newcommand{\varep}{\varepsilon}

\DeclareMathOperator{\Hom}{Hom}

\DeclareMathOperator{\SL}{SL}

\DeclareMathOperator{\ind}{ind}

\DeclareMathOperator{\disc}{disc}
\DeclareMathOperator{\height}{ht}
\DeclareMathOperator{\spn}{span}

\DeclareMathOperator{\ch}{char}
\newcommand{\triv}{\mathbf{1}}

\newcommand{\ZG}{\mathbf{G}} 
\newcommand{\ZT}{\mathbf{S}} 
\newcommand{\Tor}{S} 
\newcommand{\tor}{\mathfrak{s}} 
\newcommand{\ZX}{\mathbf{X}} 

\newcommand{\Gk}{\mathcal{G}}
\newcommand{\Vk}{\mathcal{V}}

\newcommand{\Gm}{G_{m, 0}} 
\newcommand{\Vm}{V_{m}} 
\newcommand{\ZVm}{\mathbf{V}_m} 
\newcommand{\Deltam}{\Delta_m} 

\newcommand{\Ge}{G_{8, 0}}
\newcommand{\Ve}{V_{8}} 
\newcommand{\ZVe}{\mathbf{V}_8} 
\newcommand{\Deltae}{\Delta_8} 

\newcommand{\Gs}{G_{6, 0}}
\newcommand{\Vs}{V_{6}} 
\newcommand{\ZVs}{\mathbf{V}_6} 
\newcommand{\Deltas}{\Delta_6} 
\newcommand{\ZGs}{\mathbf{G}_{6, 0}}

\newcommand{\Gf}{G_{4, 0}}
\newcommand{\Vf}{V_{4}} 
\newcommand{\ZVf}{\mathbf{V}_4} 
\newcommand{\Deltaf}{\Delta_4} 

\newcommand{\Gt}{G_{3, 0}}
\newcommand{\Vt}{V_{3}} 
\newcommand{\ZVt}{\mathbf{V}_3} 
\newcommand{\Deltat}{\Delta_3} 

\theoremstyle{plain}
\input xy
\xyoption{all}
\newtheorem{Thm}{Theorem}[section]
\newtheorem{Lem}[Thm]{Lemma}
\newtheorem{Prop}[Thm]{Proposition}
\newtheorem{corollary}[Thm]{Corollary}

\newtheorem{conjecture}[Thm]{Conjecture}
\theoremstyle{remark}
\newtheorem{Rem}[Thm]{Remark}
\numberwithin{equation}{section}
\numberwithin{paragraph}{section}

\AtEndDocument{\bigskip{\footnotesize 
\textsc{Beth Romano}~~ \texttt{beth.romano@kcl.ac.uk}\\\textsc{Department of Mathematics, King's College London,
WC2R 2LS, UK}}}

\begin{document}  

\author{Beth Romano}
\title{Stable vectors in dual Vinberg representations of $F_4$}
\maketitle

\thispagestyle{empty}  

\begin{abstract}
This paper gives a classification of stable vectors in dual Vinberg representations coming from a graded Lie algebra of type $F_4$ in a way that is independent of the field of definition. Relating these gradings to Moy--Prasad filtrations, we obtain the input for Reeder--Yu's construction of epipelagic supercuspidal representations. As a corollary, this construction gives new supercuspidal representations of $F_4(\mathbb{Q}_p)$ when $p$ is small.
\end{abstract}

\tableofcontents

\section{Introduction}

Vinberg theory, the study of graded semisimple Lie algebras, is well developed in the case when the Lie algebra in question is defined over an algebraically closed field of characteristic zero (cf. \cite{Vinberg}, \cite{Panyushev}) or of sufficiently large positive characteristic (cf. \cite{Levy}, \cite{RLYG}). In these cases we have uniform theorems describing the invariant theory of so-called Vinberg representations:
for example \cite[Theorem 3.5]{Panyushev}, which describes the ring of invariants of a Vinberg representation in terms of the invariants of the adjoint representation.
Yet when working with graded Lie algebras over fields of small characteristic, these theorems no longer hold: for example, the relevant rings of invariants need not be as described in \cite{Panyushev} (see, e.g., \cite[Section 7.5]{ReederYu}). While special cases have been studied in small characteristic  (see, e.g., \cite{Amicone}), we are still building our knowledge of this setting on a case-by-case basis. 

In this paper, I describe stable vectors (in the sense of geometric invariant theory) in certain representations coming from graded Lie algebras of type $F_4$. My results are completely general in that they describe stable vectors 
in a way that does not depend on the field of definition. These results are motivated by the construction of supercuspidal representations of $p$-adic groups in \cite{ReederYu}. For applications to $p$-adic groups, it is particularly important to classify stable vectors in the small-characteristic case, as these lead to new supercuspidal representations that do not arise from any other known construction.
I will now further describe the setting for these applications. 

Let $\Gk$ be a split, absolutely simple connected reductive group over a nonarchimedean local field $k$ with residue field $\frakf$ of characteristic $p$. The supercuspidal representations of $\Gk(k)$ have been completely classified in the case when $p$ does not divide the order of the Weyl group of $\Gk$: by \cite{Kim} and \cite{Fintzen}, all such representations may be obtained by the construction of \cite{Yu}. If we further assume that $\Gk$ is a classical group, then in the case when $p \neq 2$, all supercuspidal representations of $\Gk(k)$ may be obtained by the construction of \cite{Stevens}. Yet if $\Gk$ is an exceptional group and $p$ does divide the order of the Weyl group of $\Gk$, then there is no construction yielding all supercuspidal representations.
In \cite{ReederYu}, Reeder--Yu gave a new construction of certain \textit{epipelagic} supercuspidal representations. This construction has the benefit of being uniform for all primes $p$, but the input for their construction (which I will describe below) is not yet fully understood when $p$ is small. 

In this paper, I describe this input in the case when $\Gk$ is of type $F_4$. 
To do so, I analyze representations coming from certain graded Lie algebras. 
To describe the results in more detail, 
let $\mathcal{A}$ be the apartment associated to a maximal split torus of $\Gk$ in the Bruhat--Tits building of $\Gk$, and
let $K$ be the maximal unramified extension of $k$ with residue field $\bF$. Given a rational point $x \in \mathcal{A}$, we may form Moy--Prasad-filtration subgroups (see \cite{MP1}, \cite{MP2})
\begin{equation*}
\Gk(k)_{x, 0} \gneq \Gk(k)_{x, r_1} \gneq \Gk(k)_{x, r_2} \gneq ...
\end{equation*}
where the $r_i$ are rational numbers that depend on $x$. The quotient $\Gk(k)_{x, 0}/\Gk(k)_{x, r_1}$ forms the $\frakf$-points of a connected reductive group $\Gk_x$ over $\frakf$, and the $\frakf$-vector space $\Vk_x := \Gk(k)_{x, r_1}/\Gk(k)_{x, r_2}$ forms an algebraic representation of $\Gk_x$ with action induced by conjugation in $\Gk$. 
We say a vector $v$ in the dual representation $\check \Vk_x$ is \textit{stable} if it is stable in the sense of geometric invariant theory as a vector in $\check \Vk_x(\bF) = \check \Vk_x \otimes_\frakf \bF$ for the action of $(\Gk_x)_{\bF}$.

In \cite{ReederYu}, Reeder--Yu give a construction of supercuspidal representations of $\Gk(k)$ whose only ingredient is a stable vector in $\check \Vk_x$. To implement their construction, one must answer the following questions:
\begin{itemize}
\item[1.] For which $x$ does $\check \Vk_x$ contain stable vectors?
\item[2.] If $\check \Vk_x$ contains stable vectors, how can we describe them?
\end{itemize}
To approach these questions, we use the fact that the representation of $\Gk_x$ on $\Vk_x$ may be realized as coming from a graded Lie algebra over $\frakf$. Using Vinberg theory, 
we obtained a classification (given by \cite{ReederYu} in the case of $p$ sufficiently large, and proven uniformly for all $p$ in \cite{FintzenRomano}) of all $x$ such that $\check \Vk_x(\bF)$ contains stable vectors for the action of $(\Gk_x)_{\bF}$. 
For $p$ large enough, if $\check \Vk_x(\bF)$ contains stable vectors, then so does $\check \Vk_x$ \cite[Corollary 6.1]{Reeder}. In every known example, if $\check \Vk_x(\bF)$ contains stable vectors, then so does $\check \Vk_x$, but it is unknown whether this holds in general for arbitrary $p$. 

Similarly, the answer to the second question is well understood when $p$ is large enough: stable vectors may be described in terms of regular semisimple elements of a certain graded Lie algebra \cite[Lemma 13]{RLYG}. But this classification does not hold when $p$ is small (the vectors described in Section \ref{section-dual-cyclic} give counterexamples). There is little known about the answer to question 2 for the case when $p$ is small. 
Except for certain examples (e.g. when $x$  is the barycenter of an alcove in $\mathcal{A}$ \cite[Section 2.6]{ReederYu}; and for all $x$ when $\Gk$ is of type $G_2$ \cite{RomanoThesis}), the answer to question 2 is still unknown for arbitrary $p$. 

In this paper, I answer questions 1 and 2 in the case when $\Gk$ is of type $F_4$, for all $x$ except for one orbit under the affine Weyl group.  
To state the results, let $(X^*, \Phi, {X}_*, \check\Phi)$ be the root datum of $\Gk$, and let $\check\rho$ be the sum of a set of fundamental coweights for some choice of simple roots in $\Phi$. Choose a hyperspecial point $x_0 \in \mathcal{A}$, which determines an identification of $\mathcal{A}$ with ${X}_* \otimes \mathbb{R}$. 

\begin{Thm}
Suppose $\Gk$ is of type $F_4$ and $x$ is a rational point in $\mathcal{A}$. Assume that $x$ is not conjugate under the affine Weyl group to $x_0 + \frac{1}{2}\check\rho$.  
\begin{enumerate}
\item[1.]  If $\check{\Vk}_x(\bF)$ contains stable vectors, then so does $\check{\Vk}_x$.
\item[2.] There exists an invariant polynomial $\Delta_x$ on $\check{\Vk}_x(\bF)$ such that $v \in \check{\Vk}_x(\bF)$ is stable under the action of $(\Gk_x)_{\bF}$ if and only if $\Delta_x(v) \neq 0$.
\end{enumerate}
\end{Thm}

Moreover, we explicitly describe the invariant polynomial $\Delta_x$ in every case. 
As a corollary, the construction of \cite{ReederYu} may be used when $k = \QQ_p$ and $\Gk$ is of type $F_4$ to construction supercuspidal representations of $\Gk(\QQ_p)$ for all $p$. When $p$ is small, we obtain new supercuspdial representations.

In fact, my results are more general than stated above: for all $x$ such that $\check{\Vk}_x(\bF)$ contains stable vectors, there is a group $\mathbf{G}_x$ defined over $\bZ$ and an algebraic representation $\mathbf{V}_x$ of $\mathbf{G}_x$ over $\bZ$ such that  $(\mathbf{G}_x)_\mathfrak{f} = \Gk_x$ and $\mathbf{V}_x(\mathfrak{f}) = \check \Vk_x$. The polynomial $\Delta_x$ is defined on $\mathbf{V}_x$, and for any algebraically closed field $L$, a vector $v \in \mathbf{V}_x(L)$ is stable if and only if $\Delta_x(v) \neq 0$. This means that the description of stable vectors is uniform over any field.

We finish the introduction by describing a conjecture that would give many new examples of stable vectors in the representations $\check{\Vk}_x$. With notation as above, for any $m \in \bZ_{> 0}$, let $x_m = x_0 + \frac{1}{m}\check\rho$. (Note that if $\check{\Vk}_x$ contains stable vectors, then $\check{\Vk}_x \simeq \check{\Vk}_{x_m}$ for some $m$.) Let $h$ be the Coxeter number of $\Gk$. If $m \mid h$, then there is a natural embedding 
$\iota_m: \check{\Vk}_{x_h} \hookrightarrow \check{\Vk}_{x_m}$ (see Section \ref{section-dual-cyclic}). The stable vectors in $\check{\Vk}_{x_h}$ are well understood \cite[Section 2.6]{ReederYu}, and the Reeder--Yu construction with these vectors can be seen as the prototype for the general case (this case was first considered in the earlier paper \cite{GrossReeder}). Thus it is natural to ask about the relation between stable vectors in $\check{\Vk}_{x_h}$ and $\check{\Vk}_{x_m}$, as well as about the potential implications for representations of $p$-adic groups.

If $p$ is large, then whenever $\lam \in \check \Vk_{x_h}$ is stable for the action of $\Gk_{x_h}$, $\iota(\lam)$ is stable for the action of $\Gk_{x_m}$ \cite[Lemma 13]{RLYG}.  
If $p$ is small, then the image of a stable vector under $\iota$ may not be stable as an element of $\check \Vk_{x_m}$ \cite[Section 4.4]{RomanoThesis}. Yet all known examples satisfy the following conjecture.
\begin{conjecture}\label{conj}
Suppose $v \in \check{\Vk}_{x_h}$ is stable for the action of $\Gk_{x_h}$. Then $\iota_m(v)$ is stable for the action of $\Gk_{x_m}$ if and only if $\ch(\frakf) \nmid \frac{h}{m}$.
\end{conjecture}
If we can prove the conjecture, then, firstly, we will be able to answer question 1 in many new cases. Another reason to explore $\iota_m(v)$ for stable $v$ is because 
our understanding of Langlands parameters corresponding to simple supercuspidal representations
\cite{GrossReeder} may help to construct Langlands parameters for the supercuspidal representations associated to the vectors $\iota_m(v)$
(see \cite{RomanoThesis}). Conjecture \ref{conj} is true when $\Gk$ is of type $G_2$ \cite[Section 4.4]{RomanoThesis}. In the present paper, we verify the conjecture when $\Gk$ is of type $F_4$, in all but one case.
\begin{Thm}\label{thm-conj}
If $\Gk$ is of type $F_4$, then Conjecture \ref{conj} holds for all $m \neq 2$. 
\end{Thm}

To clarify, we restrict to the special case of $m \mid h$ in the conjecture because this is the only case in which we have a natural embedding $\check{\Vk}_{x_h} \hookrightarrow \check{\Vk}_{x_m}$. But for any $x$, stable vectors in $\check{\Vk}_x$ give the input for the construction of \cite{ReederYu} and thus have applications to the representation theory of $p$-adic groups. For $F_4$, by \cite{FintzenRomano}, we know that stable vectors exist in $\check{\Vk}_x(\bF)$ if and only if $x$ is conjugate under the affine Weyl group to $x_m$ when $m \in \{2, 3, 4, 6, 8, 12\}$. Here we consider each of these choices of $m$, except for the case when $m = 2$.

\begin{Rem}
Note that the methods described in this paper should be sufficient to extend the theorems above to the missing case, i.e. to the case when $x$ is conjugate under the affine Weyl group to $x_0 + \frac{1}{2}\check\rho$, but only after one identifies an appropriate invariant polynomial $\Delta_x$. In this case, we have not yet been able to construct such a polynomial. 
The representation $\check{\Vk}_x$ is more complicated in this case than in the other cases that appear for $F_4$, and while some of the necessary invariant theory is understood (see, e.g., \cite{Laga}), we would need an explicit analogue of the polynomial $\Delta$ in \cite[Section 3.6]{Laga} to proceed as in the other cases. 
An understanding of stable vectors in this missing case would have applications to the representation theory of $p$-adic groups in the same way as in the other cases.
\end{Rem}

\subsection*{Structure}

In Section \ref{section-one-dim}, we review the notion of stability and  a criterion for a vector to be stable. In Section \ref{section-sl-reps}, we review principal gradings on simple Lie algebras, expanding on the conjecture mentioned in the introduction and giving some useful lemmas for identifying representations. We then specialize to the group $F_4$ and its gradings. In each of sections \ref{section-3} -- \ref{section-8}, we fix $m$. We describe the representation $\check{\Vk}_{x_m}$, define an invariant polynomial $\Deltam$ on it, and then prove that  $v \in \check{\Vk}_{x_m}$ is stable if and only if $\Deltam(v) \neq 0$. Section \ref{section-padic} describes applications to $p$-adic groups, as outlined above.
For easy reference, the appendix in Section \ref{section-disc} gives formulas for certain discriminant polynomials: 
these discriminants are used to construct the invariant polynomials $\Deltam$ throughout the paper. 

\subsection*{Notation} Given $n \in \ZZ_{> 0}$ and $j \in \ZZ_{\geq 0}$, define $P_j(n)$ to be the natural representation of $\SL_n$ on homogeneous degree-$j$ polynomials in $n$ variables. If these variables are given by $X_1, ..., X_n$, we'll sometimes write $P_j(n)$ as $P_j(X_1, ..., X_n)$ for clarity. Similarly, we sometimes write $\SL_n$ as $\SL_{(X_1, \dots, X_n)}$ to specify that $\SL_n$ is acting on $P_j(X_1, \dots, X_n)$. Note that the representation $P_j(n)$ is defined over $\ZZ$. If the representation of $\SL_n$ on $P_j(n)$ factors through a quotient $\SL_n/Z_0$ for some finite central sub-group scheme $Z_0$, then we also write $P_j(n)$ for the corresponding representation of $\SL_n/Z_0$.

For $n \in \{2, 3\}$, we write $\disc_j$ for the discriminant polynomial on $P_j(n)$, which is an invariant polynomial for the action of $\SL_n$. For clarity, we sometimes write $\disc_{j, (X_1, \dots, X_n)}$ for the discriminant on $P_j(X_1, \dots, X_n)$. Note that for $f \in P_j(n)$, $\disc_j(f)$ is a polynomial over $\ZZ$ in the coefficients of $f$. 
 If $V$ is a representation of some group $G$, we write $\check{V}$ for the dual representation of $G$.

\subsection*{Acknowledgements}

The author received support from EPSRC First Grant EP/N007204/1. Many of the ideas in this paper grew out of work in my PhD thesis, and I would like to thank my advisor, Mark Reeder, for many helpful conversations related to these ideas. I'd also like to thank Jessica Fintzen for helpful conversations while collaborating on related work and Jef Laga for useful discussions. 

\section{Preliminaries}\label{section-prelim}

Throughout this section, let $G$ be a split connected reductive group of rank $\ell$ over some field $L$. Choose a split maximal torus $\Tor$ of $G$. Let $(X^* = \Hom(\Tor, \mathbb{G}_m), \Phi, X_*, \check\Phi)$ be the root datum of $G$ corresponding to $\Tor$.

\subsection{A criterion for stability}\label{section-one-dim}

In this subsection, assume $L$ is algebraically closed and that $V$ is a degree-$n$ representation of $G$ over $L$. 
Recall that a vector $v \in V$ is called \textit{stable} if its orbit $G \cdot v$ is closed and its stabilizer in $G$ is finite.

We now review a proposition, based on the Hilbert--Mumford Criterion \cite{Mumford}, that allows us to identify stable vectors. Choose a basis $\{v_1, \dots, v_n\}$ of $V$ on which $\Tor$ acts diagonally. There are $\bZ$-linear maps $\gamma_i: {X}_* \to \bZ$ for $i \in \{1, \dots, n\}$ such that for all $\lam \in {X}_*$ and $t \in L^\times$, we have $\lam(t)\cdot v_i = t^{\gamma_i(\lam)}v_i$. Given a vector $v = \sum_{i = 1}^n a_iv_i \in V$, the set $\{\gamma_i(\lam) \mid a_i \neq 0\}$ is called the set of \textit{weights} for $v$ with respect to $\lam$.
Given $\lam \in {X}_*$, let $V_\lam = \spn_L \{v_i \mid \gamma_i(\lam) \geq 0\}$, so $V_\lam$ is the set of vectors that have no negative weights with respect to $\lam$.

\begin{Prop}\label{prop-hilmum}
Let $\Delta$ be an invariant polynomial on $V$ for the action of $G$. 
Suppose $\Delta \mid_{V_\lam} = 0$ for all nontrivial $\lam \in {X}_*$. Then for any $v \in V$, if $\Delta(v) \neq 0$, the vector $v$ is stable. 
\end{Prop}
\textit{Proof}. This
follows directly from \cite[Lemma 5.1]{FintzenRomano} (see
 \cite[Remark 5.2]{FintzenRomano}). \qed

In practice, if $n$ is not small, the hypothesis of Proposition \ref{prop-hilmum} takes a long time to check. To reduce computations, we prove the following lemma.

\begin{Lem}\label{lem-weyl}
Suppose $\Delta$ is an invariant polynomial on $V$ for the action of $G$, that $\lam \in {X}_*$, and that $w$ is in the Weyl group $N_G(\Tor)/\Tor$. Then if $\Delta\mid_{V_\lam} = 0$, we have $\Delta\mid_{V_{w\lam}} = 0$. 
\end{Lem}

\textit{Proof}. Choose a lift $n \in N_G(\Tor)$ of $w$. 
Suppose the basis vector $v_j \in V_{w\lam}$. Note that $\lam(t)\cdot n^{-1}v_j = n^{-1}(w\lam)(t)n n^{-1}v_j = n^{-1}(w\lam)(t)v_j = t^{\gamma_j(w\lam)}n^{-1}v_j$, so $n^{-1}v_j \in V_{\lam}$. Thus if $v \in V_{w\lam}$, then $n^{-1}v \in V_\lam$, so by the hypothesis, $\Delta(n^{-1}v) = 0$. Since $\Delta$ is $G$-invariant, we have that $\Delta(v) = 0$. \qed

\begin{Rem}Note that Lemma \ref{lem-weyl} allows for a subtantial improvement in the amount of case-by-case analysis necessary to check the criterion of Proposition \ref{prop-hilmum}. As a comparison, see the proof of \cite[Lemma 2.16]{RomanoThesis}, which could be reduced by half if we use the lemma.
\end{Rem}

Finally, we recall a special case of the Hilbert--Mumford Criterion that will help us prove certain vectors are not stable.

\begin{Lem}\label{lem-notstable}
Let $v \in V$. If there exists a nontrivial cocharacter $\lam \in X_*$ such that $v$ has no negative weights with respect to $\lam$, then $v$ is not stable.
\end{Lem}

\textit{Proof}. This follows directly from the Hilbert--Mumford Criterion \cite[Section 1.1]{Mumford}. \qed

\subsection{On dual Vinberg representations}\label{section-sl-reps}

We return to the notation at the beginning of the section, so that $G$ is a split  connected reductive group over some field $L$ with split maximal torus $\Tor$ and $\Phi$ is the set of roots of $G$ with respect to $\Tor$.  For the rest of Section \ref{section-prelim}, we assume in addition that $G$ is absolutely simple.
 Let $\frakg$ be the Lie algebra of $G$, and let $\tor$ be the Lie algebra of $\Tor$.  Given $\al \in \Phi$, let $\frakg_\al$ be the root space corresponding to $\al$. Choose a basis $\{\al_i \mid 1 \leq i \leq \ell \}$ of simple roots of $\Phi$, and let $\check\rho \in {X}_* \otimes \QQ$ be the sum of the simple coweights with respect to this basis. Choose a $\ZZ$-form $\ZG$ of $G$ corresponding to a Chevalley basis $\{h_i, e_\al \mid 1 \leq i \leq \ell, \al \in \Phi\}$
of $\frakg$, where $\{h_i \mid 1 \leq i \leq \ell\}$ is a basis of $\text{Lie} (S)$. Let $\frakg_{\bZ} = \spn_{\bZ} \{h_i, e_\al \mid 1 \leq i \leq \ell, \al \in \Phi\}$ be the Lie algebra of $\ZG$, and let $\tor_{\bZ} = \spn_{\bZ}\{h_i \mid 1 \leq i \leq \ell\}$. Write $\mathbf{S} \subset \mathbf{G}$ for the maximal torus corresponding to $\tor_\bZ$.
Given $\al \in \Phi$, write $X_\al \subset G$ (respectively $\ZX_\al \subset \ZG$) for the corresponding root group.
Given $\al \in \Phi$, we write $G_\al \subset G$ (respectively $\ZG_\al \subset \ZG$) for the subgroup of type $A_1$ generated by the root groups $X_\al$ and $X_{-\al}$ (respectively $\ZX_\al$ and $\ZX_{-\al}$).

Given $x \in {X}_* \otimes \QQ$, we may define a grading on the Lie algebras $\frakg_{\bZ}, \frakg$ as follows. Let $m \in \bZ_{> 0}$ be minimal such that $\al(x) \in \frac{1}{m}\bZ$ for all $\al \in \Phi$. For $j \in \bZ/m\bZ$, define $(\frakg_{\bZ})_{x, j}$ as follows:
 \begin{eqnarray*}\label{eqn-xgrading}
(\frakg_\ZZ)_{x, 0} &=& \tor_\ZZ \oplus \spn_\ZZ \{e_\al \mid \al(x) \in \ZZ\}\\
(\frakg_\ZZ)_{x, j} &=& \spn_\ZZ \{ e_\al \mid \al(x) \equiv \frac{j}{m} \text{ mod } \ZZ \} ~~ (j \neq 0).
 \end{eqnarray*}
Setting $\frakg_{x, j} = (\frakg_\bZ)_{x, j} \otimes_\bZ L$ yields a $\bZ/m\bZ$-grading \begin{equation*}
\frakg = \underset{j \in \bZ/m\bZ}\bigoplus \frakg_{x, j}.
\end{equation*} 
 It's not hard to see that up to isomorphism, this grading depends only on the orbit of $x$ under the action of the affine Weyl group of $G$. Let $G_{x, 0}$ be the connected subgroup of $G$ with Lie algebra $\frakg_{x, 0}$, and define $\ZG_{x, 0} \subset \ZG$ similarly. Note that if $L = k$ is a nonarchimedean local field as in the introduction (so $G = \Gk$ as defined there), then under the identification of $\mathcal{A}$ with ${X}_* \otimes \RR$ given by a choice of hyperspecial point and with notation as in the introduction we have $\Gk_x \simeq (\ZG_{x, 0})_\frakf$ (where the subscript $\frakf$ indicates base change) and $\Vk_x \simeq (\frakg_\ZZ)_{x, 1} \otimes_\bZ \frakf$ (this follows from the proofs in \cite[Section 3.3.2]{RomanoThesis}).

\subsubsection{Dual cyclic vectors}\label{section-dual-cyclic}

Given $m \in \bZ_{> 0}$, let $x_m = \frac{1}{m}\check\rho$. Given $\al \in \Phi$,
write $\height(\al)$ for the height of $\al$ with respect to the root basis $\{\al_i \mid 1 \leq i \leq \ell\}$ (so if $\al = \sum_{i = 1}^\ell a_i\al_i$, then $\height(\al) = \sum a_i$). 
Taking $x = x_m$ in the definition above yields
\begin{eqnarray*}
\frakg_{x_m, 0} &=& \tor + \sum_{\height(\al) \equiv 0 \mod m} \frakg_\al\\
\frakg_{x_m, j} &=& \sum_{\height({\al}) \equiv j \mod m} \frakg_\al ~~(j \neq 0).
\end{eqnarray*}

Write $h$ for the Coxeter number of $G$. 
If $m \mid h$, then $\frakg_{x_h, 1} \subset \frakg_{x_m, 1}$. 
The set $\{ e_\al \mid \height(\al) \equiv 1 \mod m\}$ forms a basis of $\frakg_{x_m, 1}$. Let $\{f_\al\}$ be the dual basis of $\check\frakg_{x_m, 1}$, and denote by $f'_\al$ the restriction of $f_\al$ to $\frakg_{x_h, 1}$. Then $\{f'_\al \mid \height(\al) \equiv 1 \mod h \}$ forms a basis of $\check\frakg_{x_h, 1}$. Let us call $v \in \check\frakg_{x_m, 1}$ \textit{dual cyclic} if $v$ is $G_{x_m, 0}$-conjugate to a vector of the form
\begin{equation*}
\sum_{\height(\al) \equiv 1 \mod h} c_\al f_\al
\end{equation*}
for some constants $c_\al$ with $\prod c_\al \neq 0$ (the term \textit{cyclic} refers to the cyclic elements of \cite[Section 6.2]{Kostant}). If $L$ is algebraically closed, then $\lam = \sum c_\al f'_\al \in \check\frakg_{x_h, 1}$ is stable for the action of $G_{x_h, 0}$ if and only if $\prod c_\al \neq 0$, so a vector $v$ is dual cyclic if and only if it is $G_{x_m, 0}$-conjugate 
to the image of a stable vector under the natural embedding $\check\frakg_{x_h, 1} \to \check\frakg_{x_m, 1}$ defined by $f'_\al \mapsto f_\al$. This embedding is the one denoted $\iota_m$ in the introduction.

\subsubsection{Subrepresentations}

In this subsection, assume that $G$ (and thus $\ZG$) is adjoint. Here we'll prove some lemmas that will allow us to identify the representations of $G_{x_m, 0}$ on $\check\frakg_{x_m, 1}$. 
For the notation $P_n(k)$, see the introduction.

\begin{Lem}\label{lem-SL2-reps}
 Given $\be_0 \in \Phi$, let $\Phi_0 = \{\be_0, ..., \be_n\}$ be the root string of $\al$ through $\be_0$, where $\be_j = \be_0 + j\al$ for all $j$. Let $V = \spn_{\bZ}\{ e_\be \mid \be \in \Phi_0\}$, an algebraic representation of $\ZG_\al$ over $\ZZ$. Then  $\ZG_\al \simeq \SL_2/Z_0$ for some finite central sub-group scheme $Z_0$ of $\SL_2$, and under this isomorphism $\check{V} \simeq P_n(2)$.
\end{Lem}

\textit{Proof}. 
Note that $n \in \{0, 1, 2, 3\}$. If $n = 0$, the lemma is trivial, and if $n = 3$, the lemma follows from the proof of \cite[Proposition 4.1]{RomanoThesis}. 
Suppose $n = 1$. Given $\gamma \in \Phi$, let $u_\gamma: \mathbb{G}_a \to \ZX_\al$ be the corresponding root-group homomorphism. Note that $[e_\al, e_{\be_0}] = \varep e_{\be_1}$ for some $\varep \in \{\pm 1\}$. Then with respect to the basis $\{e_{\be_0}, e_{\be_1}\}$,  $u_\al(t)$ acts on $V$ via $\begin{psmallmatrix} 1 & 0\\ t\varep & 1 \end{psmallmatrix}$, and $u_{-\al}(t)$ acts on $V$ via $\begin{psmallmatrix} 1 & \varep t\\ 0 & 1\end{psmallmatrix}$, so we see that the representation of $\ZG_\al$ on $V$ induces an isomorphism $\ZG_\al \simeq \SL_2$. 
It's easy to check that the linear map $\phi: \check V \to P_1(X, Y)$ given by $\phi(f_{\be_0}) = X, \phi(f_{\be_1}) = -\varep Y$ is equivariant with respect to $u_\al(t)$ and $u_{-\al}(t)$, and thus $\phi$ gives the desired $\ZG_\al$-equivariant isomorphism. 
When $n = 2$, a similar straightforward calculation  gives the result. \qed

\begin{Lem}\label{lem-SL3-reps}
Suppose $\Phi_H \subset \Phi$ is a sub-root system of type $A_2$ with root basis $\{\al, \be\}$. Let 
 $H$ be the corresponding subgroup of $\ZG$ of type $A_2$ generated by root groups $X_{\pm \al}, X_{\pm \be}, X_{\pm (\al + \be)}$.
Suppose $\Phi_0 \subset \Phi$ is a subset closed under addition by $\Phi_H$, and let $V:= \spn\{ e_\del \mid \del \in \Phi_0\}$ be the corresponding subrepresentation of the representation of $H$ on $\frakg_\ZZ$. Then $H \simeq \SL_3/Z_0$ for some finite central sub-group scheme $Z_0$ of $\SL_3$, and under this isomorphism we have the following:
\begin{enumerate}
\item[1.] If $\Phi_0 = \{\gamma, \gamma + \al, \gamma + \be\}$ for some $\gamma \in \Phi$, then $\check{V} \simeq P_1(3)$.
\item[2.] If $\Phi_0 = \{\gamma, \gamma + \al, \gamma + 2\al, \gamma + \al + \be, \gamma + \al + 2\al + \be, \al + 2\al + 2\be\}$ for some $\gamma \in \Phi$, then $\check{V} \simeq P_2(3)$.
\end{enumerate}
\end{Lem}

The proof is a straightforward computation similar to the proof of Lemma \ref{lem-SL2-reps}.

\subsection{Fixing notation}\label{sec-notation}

We maintain our notation from the previous section, but from now on, we assume that $G$ is of type $F_4$ and that $L$ is algebraically closed. We make no assumption about the characteristic of $L$. We assume the simple roots $\al_1, \dots, \al_4$ are chosen to correspond to the following nodes on the Dynkin diagram of $F_4$:
$\al_1 ~\al_2 \Rightarrow \al_3 ~\al_4$. We write $\{\check\omega_1, \dots, \check\omega_4\}$ for the set of simple coweights with respect to our choice of simple roots. Given $\al \in \Phi$, we write $w_\al$ for the corresponding simple reflection in the Weyl group of $G$ with respect to $\Tor$.

To ease calculations, we make use of normalized Kac coordinates, which we write as labelings of the affine Dynkin diagram of $F_4$. For background on Kac coordinates, see \cite{ReederTA}, and to see them in the context of graded Lie algebras, see \cite{RLYG}. In each of sections \ref{section-3} -- \ref{section-8}, we use the following set-up. 
 Fix $m \in \{2, 3, 4, 6, 8\}$ ($m$ will be different in each section below). The point $\frac{1}{m}\check\rho \in {X}_* \otimes \RR$ is conjugate under the affine Weyl group of $F_4$ to an element $y_m \in {X}_* \otimes \RR$ whose Kac coordinates $a_0 ~ a_1 ~ a_2 \Rightarrow a_3 ~ a_4$ satisfy $s_i \in \{0, 1\}$ for all $i$. These Kac coordinates are listed in \cite[Table 6]{RLYG}. 
Write $\Gm$ for $G_{y_m, 0}$ as defined in Section \ref{section-sl-reps}, and similarly for $\ZG_{m, 0}$. Let $\Phi_m = \{\al \in \Phi \mid \al(x_m) \in \ZZ\}$, so the root datum of $\Gm$ is given by $(X^*, \Phi_m, {X}_*, \check{\Phi}_m)$.

We write $\ZVm$ for the representation of $\ZG_{m, 0}$ dual to $(\frakg_\ZZ)_{y_m, 1}$, and $\{ f_\al \mid \al(y_m) \equiv \frac{1}{m} \text{ mod } \ZZ \}$ for the basis dual to the natural basis $\{ e_\al \mid \al(y_m) \equiv \frac{1}{m} \text{ mod } \ZZ \}$. For any field $F$, we write $\ZVm(F)$ for the base change $\ZVm \otimes_{\ZZ} F$. 
For ease of notation, we write  $\Vm$ for the representation $\ZVm(L)$ of $\Gm$. For $q$ a power of a prime integer $p$, we say that a vector $v \in \ZVm(\bF_q)$ is \textit{stable} if it is stable as an element of $\ZVm(\overline{\bF}_p)$ for the action of the base change $(\ZG_{m, 0})_{\overline{\bF}_p}$. 

In each case below, we will construct a polynomial on $V_m$ that is invariant for the action of $\Gm$ and satisfies the hypothesis of Proposition \ref{prop-hilmum}. Note that when $m = 6$, we must first work with the representation $\ZVm$ to construct an appropriate invariant polynomial on $V_m$, but in the other cases, we work directly with $V_m$.

\section{Stable 3-grading}\label{section-3}

The normalized Kac coordinates for $\frac{1}{3}\check\rho$ are $0~0~1 \Rightarrow 0~0$, so
 $\Gt$ is of type $A_2 \times A_2$, and thus is a quotient of $\SL_3 \times \SL_3$ by a finite central subgroup. Then $\Vt$ is a representation of $\SL_3 \times \SL_3$ via precomposition with the natural quotient map, and by the definition of stability, a vector $v \in \Vt$ is stable for the action of $\Gt$ if and only if $v$ is stable for the action of $\SL_3 \times \SL_3$. 
For simplicity, we work with $\SL_3 \times \SL_3$ instead of $\Gt$ in the rest of the section. 
Lemma \ref{lem-SL3-reps} implies that $\Vt \simeq P_1(3) \otimes P_2(3)$ as a representation of $\SL_3 \times \SL_3$. 

\subsection{An invariant polynomial}

We think of vectors $v \in \Vt$ as being of the form:
\begin{equation}\label{eqn-3-v}
\begin{split}
v &= (aX + bY + cZ)\otimes T^2 + (dX + eY + fZ)\otimes U^2 + (gX + hY + iZ)\otimes W^2\\ &   + (jX + kY + lZ)\otimes TU + (mX + nY + oZ) \otimes TW +  (pX + qY + rZ)\otimes UW.
\end{split}
\end{equation}
We construct an invariant polynomial on $\Vt$  by first thinking of $v$ as a polynomial in $T, U, W$ with coefficients in $P_1(3)$. We take the discriminant $\disc_{2, (T, U, W)}$ of this polynomial, which gives a map $P_1(3) \otimes P_2(3) \to P_3(3) = P_3(X, Y, Z)$. We then take the discriminant polynomial $\disc_{3, (X, Y, Z)}$ on $P_3(3)$, yielding the composite discriminant $\Deltat := \disc_{3, (X, Y, Z)} \circ \disc_{2, (T, U, W)}$. Since $\disc_{2, (T, U, W)}: P_1(3) \otimes P_2(3) \to P_3(3)$ is invariant for the action of $1 \times \SL_3 = \SL_{(T, U, W)}$ and equivariant for the action of $\SL_3 \times 1 = \SL_{(X, Y, Z)}$, and $\disc_{3, (X, Y, Z)}$ is invariant for the action of $\SL_3$ on $P_3(3)$, we see that $\Deltat$ is an invariant polynomial for $\Vt$.

More explicitly, if $v$ is as in (\ref{eqn-3-v}), then plugging the coefficients $aX + bY + cZ, \dots, pX + qY + rZ$ into the formula (\ref{eqn-disc-tern-quad}), we see that we may write $\disc_{2, (T, U, W)}(v)$ as
\begin{equation}\label{eqn-disc2v}
AX^3 + BX^2 Y + CX^2 Z + DXY^2 + EXYZ + FXZ^2 + GY^3 + 
 HY^2 Z + IYZ^2 + JZ^3,
 \end{equation}
 where $A, B, C, ..., J$ are the following homogeneous polynomials in $a, b, c, ..., r$, invariant under the action of $\SL_{(T, U, W)}$:
 \begin{equation}\label{eqn-3A-J}
\begin{split}
A =&~ 4 a d g - g j^2 - d m^2 + j m p - a p^2\\
B =&~ 4 b d g + 4 a e g + 4 a d h - h j^2 - 2 g j k - e m^2 - 2 d m n + 
 k m p + j n p - b p^2 + j m q - 2 a p q\\
C =&~ 4 c d g + 4 a f g + 4 a d i - i j^2 - 2 g j l - f m^2 - 2 d m o + 
 l m p + j o p - c p^2 + j m r - 2 a p r\\
D =&~ 4 b e g + 4 b d h + 4 a e h - 2 h j k - g k^2 - 2 e m n - d n^2 + 
 k n p + k m q + j n q - 2 b p q - a q^2\\
 E =&~ 
4 c e g + 4 b f g + 4 c d h + 4 a f h + 4 b d i + 4 a e i - 2 i j k - 
 2 h j l - 2 g k l - 2 f m n - 2 e m o\\& - 2 d n o + l n p + k o p + 
 l m q + j o q - 2 c p q + k m r + j n r - 2 b p r - 2 a q r\\
 F =&~ 4 c f g + 4 c d i + 4 a f i - 2 i j l - g l^2 - 2 f m o - d o^2 + 
 l o p + l m r + j o r - 2 c p r - a r^2\\
 G =&~  4 b e h - h k^2 - e n^2 + k n q - b q^2\\
 H =&~ 4 c e h + 4 b f h + 4 b e i - i k^2 - 2 h k l - f n^2 - 2 e n o + 
 l n q + k o q - c q^2 + k n r - 2 b q r\\
 I =&~ 4 c f h + 4 c e i + 4 b f i - 2 i k l - h l^2 - 2 f n o - e o^2 + 
 l o q + l n r + k o r - 2 c q r - b r^2\\
 J =&~ 4 c f i - i l^2 - f o^2 + l o r - c r^2.
 \end{split}
\end{equation}
To obtain $\Deltat$, we apply $\disc_3$ to (\ref{eqn-disc2v}), using the formula obtained in Section \ref{section-disc-3-3}. The following lemma is easy to verify and will be used repeatedly in the next proof.

\begin{Lem}\label{3-disc-0}
Let
\begin{equation*}
f(X, Y, Z) = AX^3 + BX^2 Y + CX^2 Z + DXY^2 + EXYZ + FXZ^2 + GY^3 + 
 HY^2 Z + IYZ^2 + JZ^3
 \end{equation*}
 and suppose one of the following holds:
 \begin{enumerate}
 \item[1.] $D =G = H = 0$
 \item[2.] $G = H = I = J = 0$.
 \end{enumerate}
 Then $\disc_3(f) = 0$.
 \end{Lem}

The next lemma will let us show that if $\Deltat(v) \neq 0$, then $v$ is stable.
 
\begin{Lem}\label{lem-invar3}
The invariant polynomial $\Deltat = \disc_3 \circ \disc_2$ satisfies the hypothesis of Proposition \ref{prop-hilmum}.
\end{Lem}

\textit{Proof.}
Let $\lam$ be a nontrivial cocharacter of the diagonal torus of $\SL_3 \times \SL_3$. Then $\lam$ is of the form
\begin{eqnarray*}
\lam(t) = \begin{pmatrix} t^{s_1} & 0 & 0\\
0 & t^{s_2} & 0\\
0 & 0 & t^{-s_1 - s_2}
\end{pmatrix}
\times
\begin{pmatrix} t^{s_3} & 0 & 0\\
0 & t^{s_4} & 0\\
0 & 0 & t^{-s_3 - s_4}
\end{pmatrix}.
\end{eqnarray*}
for some integers $s_1, s_2, s_3, s_4$. Note that $\lam$ acts with the following weights on $\Vt$, corresponding to the variables $a, b, c, \dots, r$ in (\ref{eqn-3-v}):
\begin{eqnarray*}
\gamma_a(\lam) &=& s_1 + 2s_3\\
\gamma_b(\lam) &=& s_2 + 2s_3\\
\gamma_c(\lam) &=& -s_1 -s_2 + 2s_3\\
\gamma_d(\lam) &=& s_1 + 2s_4\\
\gamma_e(\lam) &=& s_2 + 2s_4\\
\gamma_f(\lam) &=& -s_1 -s_2 + 2s_4\\
\gamma_g(\lam) &=& s_1  -2s_3 - 2s_4\\
\gamma_h(\lam) &=& s_2 -2s_3 - 2s_4\\
\gamma_i(\lam) &=& -s_1 -s_2 -2s_3 - 2s_4\\
\gamma_j(\lam) &=& s_1 + s_3 + s_4\\
\gamma_k(\lam) &=& s_2 + s_3 + s_4\\
\gamma_l(\lam) &=& -s_1 -s_2 + s_3 + s_4\\
\gamma_m(\lam) &=& s_1 -s_4\\
\gamma_n(\lam) &=& s_2 - s_4\\
\gamma_o(\lam) &=& -s_1 -s_2 - s_4\\
\gamma_p(\lam) &=& s_1 - s_3\\
\gamma_q(\lam) &=& s_2 - s_3\\
\gamma_r(\lam) &=& -s_1 -s_2 - s_3.
\end{eqnarray*}
By this we mean $\lam(t)\cdot X\otimes T^2 = t^{\gamma_a(\lam)}X\otimes T^2$ and similarly for the other basis vectors labeled as in (\ref{eqn-3-v}).

Let $v \in (\Vt)_\lam$. Define $a, b, \dots, r$ to be the coefficients of $v$ as in (\ref{eqn-3-v}), and define $A, B, \dots, J$ to be the coefficients of $\disc_2(v)$ as in (\ref{eqn-disc2v}). We must show that $\Deltat(v) = 0$.
First suppose $s_1 = s_2 = 0$. By Lemma \ref{lem-weyl}, we may replace $\lam$ with a Weyl-conjugate, and so since $\lam$ is nontrivial we may assume $s_3 > 0$ and $s_4 < 0$. Then $\gamma_d(\lam), \gamma_e(\lam), \gamma_f(\lam), \gamma_p(\lam), \gamma_q(\lam)$, and $\gamma_r(\lam)$ are all negative. Since $v \in (\Vt)_\lam$, we have $d = e = f = p = q = r = 0$, 
so
 \begin{eqnarray*}
A &=&  g j^2  \\
B &=&   - h j^2 - 2 g j k\\
C &=&  - i j^2 - 2 g j l\\
D &=&  - 2 h j k - g k^2\\
 E &=&  - 2 i j k - 
 2 h j l - 2 g k l\\
 F &=& - 2 i j l - g l^2 \\
 G &=&  - h k^2 \\
 H &=&  - i k^2 - 2 h k l \\
 I &=&  - 2 i k l - h l^2\\
 J &=&  - i l^2.
 \end{eqnarray*}
Using the formulas in Section \ref{section-disc}, we see that $\Deltat(v) = 0$. 
 Similarly, if $s_3 = s_4 = 0$, then acting by the Weyl group, we may assume $s_1 > 0$ and $s_2 < 0$. By a similar logic to the previous case, we have $b = e = h = k = n = q = 0$, so $B = D = E = G = H = I = 0$, i.e. $\disc_2(v)$ is a polynomial in $X$ and $Z$.  Thus $\Deltat(v) = 0$.

Otherwise at least one of $s_1, s_2$ is nonzero and at least one of $s_3, s_4$ is nonzero. Acting by the Weyl group of $\SL_3 \times \SL_3$, we may assume $s_1, s_3 > 0$ and $s_2, s_4 < 0$. 
This implies $\gamma_e(\lam) < 0$ and $\gamma_q(\lam) < 0$. Since $v \in (\Vt)_\lam$, we have $e = q = 0$. We next claim that
\begin{itemize}
\item[C1.] $fh = 0$, 
\item[C2.] $fn = 0$,
\item[C3.] $hk = 0$, 
\item[C4.] $kr = 0$.
\end{itemize}
To prove C1, we need to show that either $\gamma_f(\lam) < 0$ or $\gamma_h(\lam) < 0$. 
But this follows from the fact that $\gamma_f(\lam) + \gamma_h(\lam) = -s_1 -2s_3 < 0$. Similarly, C2 follows from the fact that if $\gamma_n(\lam) \geq 0$, then $\gamma_f(\lam) = -s_1 + s_4 -\gamma_n(\lam) < 0$, so at least one of $\gamma_f(\lam), \gamma_n(\lam)$ is negative. For C3, note that
if $\gamma_k(\lam) \geq 0$, then $s_3 + s_4 \geq -s_2$, so $\gamma_h(\lam) = s_2 - 2s_3 - 2s_4 \leq 2s_2 < 0$. 
 Similarly, $\gamma_k(\lam) + \gamma_r(\lam) = -s_1 + s_4 < 0$, so $kr = 0$.
Thus
\begin{eqnarray*}
D &=& 4 b d h - g k^2  - d n^2 + 
 k n p\\
G &=&  0\\
H &=&   - i k^2    \\
 I &=&  4 b f i - 2 i k l - h l^2   + l n r   - b r^2\\
  J &=& 4 c f i - i l^2 - f o^2 + l o r - c r^2.
\end{eqnarray*}

We now split into three cases.

Case 1: $\gamma_k(\lam) \geq 0$.

By the same logic as in the proof of C3 above, we have $\gamma_h(\lam) < 0$, so $h = 0$. We also have $\gamma_i(\lam) \leq -s_1 + s_2 < 0$, so $i = 0$. We have $s_4 \geq -s_2 -s_3$, so $\gamma_r(\lam) < 0$ and $r = 0$. 
By the formulas above, we have $G = H = I = 0$. 

If $\gamma_g(\lam) \geq 0$, then $\gamma_g(\lam) + 2\gamma_k(\lam) = s_1 + 2s_2 \geq 0$, so $s_1 + s_2 > 0$, and thus $\gamma_f(\lam) < 0$. We see that $J = 0$, so by Lemma \ref{3-disc-0}, we have $\Deltat(v) = 0$. 

Otherwise $\gamma_g(\lam) < 0$, and we have
\begin{eqnarray*}
D &=&    - d n^2 + 
 k n p\\
  J &=&  - f o^2.
\end{eqnarray*}
If $\gamma_n(\lam) < 0$, then $D = 0$ and $\Deltat(v) = 0$ by Lemma \ref{3-disc-0}. If $\gamma_n(\lam) \geq 0$, then $\gamma_f(\lam) = -s_1 + s_4 - \gamma_n(\lam) < 0$, yielding $J = 0$ again and completing Case 1.

Case 2: $\gamma_k(\lam) < 0$ and $\gamma_n(\lam) \geq 0$.

Since $\gamma_k(\lam) < 0$ and $v \in (\Vt)_\lam$, we have $k = 0$. Since $\gamma_n(\lam) \geq 0$, by the same logic as explained above in the proof of C2, we have $f = 0$, so
\begin{eqnarray*}
D &=& 4 b d h   - d n^2\\
H &=&  0\\
 I &=&    - h l^2   + l n r   - b r^2\\
  J &=& - i l^2  + l o r - c r^2.
\end{eqnarray*}

Suppose $\gamma_l(\lam) < 0$. Then if $\gamma_r(\lam) < 0$, we have $G = H = I = J = 0$ and $\Deltat(v) = 0$ by Lemma \ref{3-disc-0}. Otherwise $\gamma_r(\lam) \geq 0$. Since $\gamma_l(\lam) < 0$ and $\gamma_n(\lam) \geq 0$, we have $s_3 < s_1$, so $0 \leq \gamma_r(\lam) = -s_1 - s_2 -s_3 < -s_2 - 2s_3$, so $\gamma_b(\lam) < 0$, which gives
\begin{eqnarray*}
D &=&    - d n^2\\
 I &=&   0  \\
  J &=&  - c r^2.
\end{eqnarray*}
Now $\gamma_d(\lam) + \gamma_c(\lam) < -\gamma_n(\lam) -\gamma_r(\lam) \leq 0$, so at least one of $c, d$ is zero, yielding either $D = 0$ or $J = 0$, so $\Deltat(v) = 0$ by Lemma \ref{3-disc-0}. 
So we may assume that $\gamma_l(\lam) \geq 0$ for the rest of Case 2.

Now $\gamma_l(\lam) \geq 0$ and $\gamma_n(\lam) \geq 0$ together imply $\gamma_d(\lam) \leq s_2 + s_3 + s_4 =  \gamma_k(\lam) < 0$, so $d = 0$. So $D = G = H = 0$, $\Deltat(v) = 0$, and we are done with Case 2. 

Case 3: $\gamma_k(\lam) < 0$ and $\gamma_n(\lam) < 0$.

Since $\gamma_n(\lam) < 0$, we have $n = 0$, so
\begin{eqnarray*}
D &=& 4 b d h \\
H &=&  0\\
 I &=&  4 b f i  - h l^2      - b r^2\\
  J &=& 4 c f i - i l^2 - f o^2 + l o r - c r^2.
\end{eqnarray*}
If $bdh = 0$, then $D = 0$ and $\Deltat(v) = 0$ by Lemma \ref{3-disc-0}, so we may assume for the rest of the proof that $\gamma_b(\lam) \geq 0, \gamma_d(\lam) \geq 0$, and $\gamma_h(\lam) \geq 0$. 
These assumptions imply that $\gamma_b(\lam) + \gamma_d(\lam) + \gamma_h(\lam) \geq 0$ and $-\gamma_n(\lam) - \gamma_k(\lam) + \gamma_b(\lam) > 0$. Adding these inequalities, we have 
\begin{equation*}
0 < 2\gamma_b(\lam) + \gamma_d(\lam) + \gamma_h(\lam) - \gamma_k(\lam) - \gamma_n(\lam) = -\gamma_r(\lam),
\end{equation*}
 so $\gamma_r(\lam) < 0$ and $r = 0$. We also have $\gamma_l(\lam) = \gamma_k(\lam) - (s_1 + 2s_2) < 0$, so $l = 0$; and $\gamma_f(\lam) = \gamma_r(\lam) + (s_3 + 2s_4) < 0$, so $f = 0$. Putting these together, we have $G = H = I = J = 0$, so $\Deltat(v) = 0$. \qed

\subsection{Stable vectors}

\begin{Thm}\label{thm-3}
A vector $v \in P_1(3) \otimes P_2(3)$ is stable for the action of $\SL_3 \times \SL_3$ if and only if $\Deltat(v) \neq 0$.
\end{Thm}

 \textit{Proof}. By Lemma \ref{lem-invar3} and Proposition \ref{prop-hilmum}, if $\Deltat(v) \neq 0$, then $v$ is stable. For the other direction, suppose $\Deltat(v) = 0$. Write $v$ with coefficients $a, b, \dots, r$ as in (\ref{eqn-3-v}), and define $A, B, \dots, J$ as in (\ref{eqn-3A-J}). Notice that if we write $v$ as
 \begin{equation*}
 v = X\otimes f_X(T, U, W) + Y\otimes f_Y(T, U, W) + Z \otimes f_Z(T, U, W),
 \end{equation*}
 then $f_X = aT^2 + dU^2 + gW^2 + jTU +  mTW + pUW$, and $A = \disc_{2, (T, U, W)}(f_X)$.
 
Since $\Deltat(v) = 0$, as in \cite[Section 1.11]{Mumford}, after acting by $\SL_{X, Y, Z}$, we may assume that $A = B = C = D = 0$. Thus $\disc_{2}(f_X) = 0$, and by \cite[Section 1.10]{Mumford}, acting by $\SL_{T, U, W}$ we may assume $a = j = m = 0$, so that $f_X = dU^2 + pUW + gW^2$.
 Plugging in, we have
\begin{eqnarray*}
B &=& b(4dg - p^2) = 0\\
C &=& c(4dg - p^2) = 0.
\end{eqnarray*}
If $4dg - p^2 \neq 0$, then $b = c = 0$, and with notation as in the proof of Lemma \ref{lem-invar3}, we may take
$(s_1, s_2, s_3, s_4) = (-2, 1, -2, 1)$ to obtain a cocharacter $\lam$ for which $v$ has no negative weights.
By Lemma \ref{lem-notstable}, $v$ is not stable. So we may assume $4dg - p^2 = 0$, and so $\disc_{2, (U, W)}(dU^2 + pUW + gW^2) = 0$ and the quadratic $dU^2 + pUW + gW^2$ is a square. 
Acting by the subgroup $\SL_{(U, W)} \subset \SL_{(T, U, W)}$, we may assume $dU^2 + pUW + gW^2 = U^2$, so $p = g = 0$. Note that moving to another element in the orbit $\SL_{(U, W)}\cdot v$ does not change the fact that $a = j = m = 0$. 
Then taking the cocharacter corresponding to $(s_1, s_2, s_3, s_4) = (-4, 2, -1, 2)$ in Lemma \ref{lem-notstable} shows that $v$ is not stable.\qed

\begin{corollary}\label{cor-3}
\begin{enumerate}
\item[1.] For every prime integer $\mathsf{p}$, the representation $\ZVt(\bF_{\mathsf{p}})$ contains stable vectors.
\item[2.] Dual cyclic vectors in $\ZVt(\bF_{\mathsf{p}})$ are stable if and only if $\mathsf{p} \neq 2$.
\end{enumerate}
\end{corollary}

\textit{Proof}.
Dual cyclic vectors are of the form $aX\otimes T^2 + eY \otimes U^2 + hY \otimes W^2 + lZ \otimes TU + pX \otimes UW$.
We have $\Deltat(aX\otimes T^2 + eY \otimes U^2 + hY \otimes W^2 + lZ \otimes TU + pX \otimes UW) = 4096 a^6 e^3 h^9 l^{12}p^6$, so if $a, e, h, l, p \in \bF_{\mathsf{p}}^\times$, then $\Deltat(aX\otimes T^2 + eY \otimes U^2 + hY \otimes W^2 + lZ \otimes TU + pX \otimes UW) \neq 0$ if and only if $\mathsf{p} \neq 2$.
This proves part 2 of the corollary. This also shows that $\Deltat$ is nonzero on $\ZVt(\FF_\mathsf{p})$ for all $\mathsf{p} \neq 2$. The fact that  
\begin{equation*}
\Deltat( (X + Z)T^2 + (X + Z) U^2 + ZW^2 + (Y + Z)TU + (X + Z)TW + YUW) = 225
\end{equation*}
shows that $\Deltat$ is nonzero on $\ZVt(\FF_2)$, proving part 1 of the corollary. 
\qed

\section{Stable 4-grading}\label{section-4}

The normalized Kac coordinates for $\frac{1}{4}\check\rho$ are $1~0~1\Rightarrow 0~0$, so we see that $\Phi_4$ is the sub-root system generated by $\{\al_1, \al_3, \al_4\}$. The subgroup of $G$ generated by $G_{\al_3}$ and $G_{\al_4}$ is of type $A_2$, and since this subgroup has an irreducible representation of dimension 6 (given by $\spn \{e_{\al_2}, e_{\al_2 + \al_3}, e_{\al_2 + 2\al_3}, e_{\al_2 + \al_3 + \al_4},\\ e_{\al_2 + 2\al_3 + \al_4}, e_{\al_2 + 2\al_3 +2\al_4}\}$), we see that $\la G_{\al_3}, G_{\al_4}\ra \simeq \SL_3$. 
Since $\al_1$ has a root string of length 2 through $\al_2$, we have $G_{\al_1} \simeq \SL_2$. 

The group $\Gf$ is generated by $G_{\al_1}, G_{\al_3}, G_{\al_4}$, and the subtorus of $\Tor$ with cocharacter group generated by $\check\omega_2$.
By the above considerations, the subgroup of $\Gf$ generated by $G_{\al_1}$, $G_{\al_3}$, and  $G_{\al_4}$ is isomorphic to $\SL_2 \times \SL_3$, and as a representation of this subgroup, we have
\begin{equation*}
\Vf \simeq (P_{1}(2) \otimes P_2(3) ) \oplus (P_{1}(2) \otimes \triv),
\end{equation*}
where $P_1(2) \otimes P_2(3)$ is generated as a subrepresentation by $f_{\al_2}$, $P_1(2) \otimes \triv$ is generated by $f_{-\al_0}$, and $\triv$ is the trivial representation of $\SL_3$.
Note that $\check\omega_2(t)$ acts on $P_1(2) \otimes P_2(3)$ as multiplication by $t^{-1}$ and acts on $P_1(2) \otimes \triv$ as multiplication by $t^3$.

\subsection{An invariant polynomial}

Choose a basis $T$ of the trivial representation of $\SL_3$. Let $v \in \Vf$ be given by $v = F_v(X, Y, Z) + (mU + nW)T$, where
\begin{equation}\label{eqn-4v}
\begin{split}
F_v(X, Y, Z) =& ~(aU + bW) \otimes X^2 + (cU + dW) \otimes Y^2 + (eU + fW) \otimes Z^2 + (gU + hW) \otimes XY\\& + (iU + jW) \otimes XZ + (kU + lW) \otimes YZ.
\end{split}
\end{equation}
The map $\Vf \to P_4(U, W)$ given by $v \mapsto (mU + nW)\disc_{2, (X, Y, Z)}(F_v)$ is equivariant for the action of $\SL_2$ and invariant under the actions of both $\SL_3$ and $\check\omega_2(t)$ for all $t \in L^\times$. 
Let 
\begin{equation*}
\Deltaf(v) = \disc_{4, (U, W)}((mU + nW)\disc_{2, (X, Y, Z)}(F_v)).
\end{equation*}
 Then $\Deltaf$ is an invariant polynomial for the action of $\Gf$ on $\Vf$. 
More explicitly, if we let 
\begin{equation*}
(mU + nW)\disc_{2, (X, Y, Z)}(F_v) = AU^4 + BU^3W + CU^2W^2 + DUW^3 + EW^4,
\end{equation*}
then
\begin{equation}\label{eqn-4A-E}
\begin{split}
A =&~ 4 a c e m - e g^2 m - a k^2 m + g i  k m- c i^2m \\
B =&~ 4 b c e m + 4 a d e m + 4 a c f m - f g^2 m - 2 e g h m - b k^2 m - 
 2 a k l m + 4 a c e n - e g^2 n - a k^2 n\\& + h i k m + g  il m + 
 g ik n  - d i^2m  - ci^2 n  + g jk m  - 2 c i jm \\
C =&~ 4 b d e m + 4 b c f m + 4 a d f m - 2 f g h m - e h^2 m - 2 b k l m - 
 a l^2 m + 4 b c e n + 4 a d e n + 4 a c f n\\& - f g^2 n - 2 e g h n - 
 b k^2 n - 2 a k l n + h il m  + h ik n + g il n  - d i^2n  + 
 h jk m  + g jl m \\& + g jk n  - 2 d i jm  - 2 ci j  n - c j^2m \\
D =&~ 4 b d f m - f h^2 m - b l^2 m + 4 b d e n + 4 b c f n + 4 a d f n - 
 2 f g h n - e h^2 n - 2 b k l n - a l^2 n\\& + h il n  + h jl m  + 
 h jk n  + g jl n  - 2 d i jn  - d j^2m  - cj^2 n \\
E =&~ 4 b d f n - f h^2 n - b l^2 n + h jl n  - dj^2 n .
\end{split}
\end{equation}
We can then find $\Deltaf$ using (\ref{eqn-quarticdisc}).

\begin{Lem}\label{lem-invar4}
The invariant polynomial $\Deltaf$ satisfies the hypothesis of Proposition \ref{prop-hilmum}.
\end{Lem}

\textit{Proof}. The proof is similar to that of Lemma \ref{lem-invar3}. As above, let $v = F_v(X, Y, Z) + (mU + nW)T$, with $F_v(X, Y, Z)$ definited as in (\ref{eqn-4v}), and define $A, B, C, D, E$ as in (\ref{eqn-4A-E}). Throughout the proof, we will use the fact that if $D = E = 0$, then $\Deltaf(v) = 0$. 

If $\lam = 
s_1\check\omega_1 + s_2\check\omega_2 + s_3\check\omega_3 + s_4\check\omega_4$
is an arbitrary cocharacter, the weights of $\lam$ corresponding to $a, b, ..., n$ are
\begin{eqnarray*}
\gamma_a(\lam) &=& -s_2\\
\gamma_b(\lam) &=& -s_1 -s_2\\
\gamma_c(\lam) &=& -s_2 - 2s_3\\
\gamma_d(\lam) &=& -s_1 -s_2 - 2s_3\\
\gamma_e(\lam) &=& -s_2 - 2s_3 - 2s_4\\
\gamma_f(\lam) &=& -s_1 -s_2 - 2s_3 - 2s_4\\
\gamma_g(\lam) &=&  -s_2 - s_3\\
\gamma_h(\lam) &=& -s_1 -s_2 - s_3\\
\gamma_i(\lam) &=& -s_2 - s_3 - s_4\\
\gamma_j(\lam) &=& -s_1 -s_2 - s_3 - s_4\\
\gamma_k(\lam) &=& -s_2 - 2s_3 - s_4\\
\gamma_l(\lam) &=& -s_1 -s_2 - 2s_3 - s_4\\
\gamma_m(\lam) &=& 2s_1 + 3s_2 + 4s_3 + 2s_4\\
\gamma_n(\lam) &=& s_1 + 3s_2 + 4s_3 + 2s_4.
\end{eqnarray*}
As in the proof of Lemma \ref{lem-invar3}, by this we mean that $\lam(t)\cdot U\otimes X^2 = t^{\gamma_a(\lam)} U \otimes X^2$ and similarly for the other basis elements with coefficients as in (\ref{eqn-4v}). Suppose $v \in (\Vf)_\lam$, i.e. that $v$ has no negative weights with respect to $\lam$. 
We claim that $\Deltaf(v) = 0$. 

By Lemma \ref{lem-weyl}, we may replace $\lam$ with an element in its orbit under the Weyl group of $\Gf$, so we may assume that $\lam$ is in the closure of the fundamental Weyl chamber in ${X}_* \otimes \RR$ with respect to the root basis $\{\al_1, \al_3, \al_4\}$ of $\Phi_4$. This means we may assume that $s_i \geq 0$ for $i \in \{1, 3, 4\}$.

Note that if $\gamma_m(\lam) < 0$, then $m = 0$, and since $s_1 \geq 0$, we have $\gamma_n(\lam) < 0$. Since $v \in (\Vf)_\lam$, we have $m = n = 0$, so $(mU + nW)\disc_{2, (X, Y, Z)}(F_v)$ is the zero polynomial and clearly $\Deltaf(v) = 0$. So for the rest of the proof we may assume $\gamma_m(\lam) \geq 0$.

If $\gamma_b(\lam) < 0$, then $b = 0$, and since $s_3 \geq 0$, we have $\gamma_d(\lam) < 0$. Since $v \in (\Vf)_\lam$, this implies $d = 0$. Similarly, we have $f = h = j = l = 0$. Plugging into  (\ref{eqn-4A-E}), we see that $D = E = 0$, and so $\Deltaf(v) = 0$. For the rest of the proof we may assume $\gamma_b(\lam) \geq 0$. 

Similarly, note that if $\gamma_c(\lam) < 0$, then $c = 0$ and $\gamma_d(\lam) < 0$, so $d = 0$. Similarly, $e = f = k = l = 0$. Plugging into  (\ref{eqn-4A-E}), we see that $D = E = 0$, and so again $\Deltaf(v) = 0$. So for the rest of the proof we may assume $\gamma_c(\lam) \geq 0$.

Next we claim that if $s_2 \geq 0$, then $\Deltaf(v) = 0$. Indeed
 if $s_2 \geq 0$, then $s_i \geq 0$ for all $i$.
Since $\lam$ is nontrivial, we must have $s_i > 0$ for some $i$. Our assumptions then imply that $\gamma_f(\lam) < 0$, and since $v \in (\Vf)_\lam$, this implies $f = 0$. Similarly, we must have $j = l = 0$. Since $\gamma_d(\lam) + \gamma_e(\lam) < 0$, we must have $de = 0$, and similarly $eh = hi = 0$. Plugging into (\ref{eqn-4A-E}), we see that $D = E = 0$, so $\Deltaf(v) = 0$. 
For the rest of the proof we assume $s_2 < 0$. 

For easy reference, we list the assumptions we have just made.
\begin{itemize}
\item[A1.] $s_1 \geq 0$.
\item[A2.] $s_2 < 0$.
\item[A3.] $s_3 \geq 0$.
\item[A4.] $s_4 \geq 0$.
\item[A5.] $\gamma_b(\lam) \geq 0$.
\item[A6.] $\gamma_c(\lam) \geq 0$.
\item[A7.] $\gamma_m(\lam) \geq 0$.
\end{itemize}

By A7, $2s_1 + 3s_2 + 4s_3 + 2s_4 \geq 0$, so $2(s_1 + 2s_3 + s_4) \geq -3s_2 > 0$, which implies $s_1 + 2s_3 + s_4 > -s_2$. Thus $s_1 + s_2 + 2s_3 + s_4 > 0$, so $\gamma_\ell(\lam) < 0$ and $\gamma_f(\lam) < 0$. Plugging in $f = l = 0$, we get
\begin{eqnarray*}
A &=& 4 a c e m - e g^2 m - a k^2 m + g k m i - c m i^2\\
B &=& 4 b c e m + 4 a d e m - 2 e g h m - b k^2 m + 4 a c e n - e g^2 n - 
 a k^2 n + h k m i + g k n i - d m i^2  - c n i^2\\& & + g k m j - 2 c m ij\\
C &=& 4 b d e m  - e h^2 m   + 4 b c e n + 4 a d e n   - 2 e g h n - 
 b k^2 n   + h ik n   - d i^2n  + 
 h jk m    + g jk n  - 2 d i jm\\& &  - 2 ci j  n - c j^2m\\
D &=& 4 b d e n - e h^2 n + h k n j - 2 d n i j - d m j^2 - c n j^2\\
E &=& -d n j^2.
\end{eqnarray*}

We now consider several special cases.

Case 1: $s_4 = 0$.

In this case $\gamma_d(\lam) = \gamma_l(\lam) < 0$. By A3 and A6, we have $-s_2 - s_3 \geq 0$, so $\gamma_n(\lam) = -\gamma_d(\lam) + 2(s_2 + s_3) < 0$. Plugging in $d = n = 0$, we see that $D = E = 0$, and $\Deltaf(v) = 0$.

Case 2: $\gamma_b(\lam) = 0$.

By Case 1 and A4 we may assume $s_4 > 0$. In this case $\gamma_j(\lam) = -s_3 -s_4 < 0$. Plugging in $j = 0$, we have $E = 0$ and $D = 4 b d e n - e h^2 n$. 
If $e = 0$, then clearly $D = 0$ and $\Deltaf(v) = 0$. If $e \neq 0$, then $\gamma_e(\lam) \geq 0$. This implies $\gamma_n(\lam) = -2\gamma_e(\lam) - 2s_4 < 0$, so $n = 0$, which implies $D = 0$ and $\Deltaf(v) = 0$.

Case 3: $\gamma_h(\lam) \leq 0$.

By Case 1 and A4, we may assume $s_4 > 0$, so $\gamma_j(\lam) = \gamma_h(\lam) - s_4 < 0$, giving $j = 0$. 
Similarly, by Case 2 and A5, we may assume $\gamma_b(\lam) > 0$, i.e. $s_1 + s_2 < 0$. Since $\gamma_h(\lam) \leq 0$, this implies $s_3 > 0$ and $\gamma_d(\lam) = \gamma_h(\lam) - s_3 < 0$. Plugging in $d = h = j = 0$, we have $D = E = 0$ and $\Deltaf(v) = 0$. This completes Case 3.

To finish, by Case 3, we may assume $\gamma_h(\lam) > 0$. Since $\gamma_m \geq 0$, this implies $s_2 + 2s_3 + 2s_4 > 0$, i.e. $\gamma_e(\lam) < 0$. Plugging in $e = 0$, we have
\begin{eqnarray*}
A &=&  - a k^2 m + g k m i - c m i^2\\
B &=&  - b k^2 m  - a k^2 n + h k m i + g k n i - d m i^2 - c n i^2 + g k m j - 2 c m ij\\
C &=&  -  b k^2 n   + h ik n   - d i^2n  + 
 h jk m    + g jk n  - 2 d i jm  - 2 ci j  n - c j^2m\\
D &=& h k n j - 2 d n i j - d m j^2 - c n j^2\\
E &=& -d n j^2.
\end{eqnarray*}
If $j = 0$, then $D = E = 0$, so $\Deltaf(v) = 0$. Otherwise $j \neq 0$, so $\gamma_j(\lam) \geq 0$, i.e. $-s_1 -s_2 -s_3 -s_4 \geq 0$. Since $\gamma_m(\lam) \geq 0$, we have $s_2 + 2s_3 \geq 0$, i.e. $\gamma_c(\lam) \leq 0$. By A6, this implies $\gamma_c(\lam) = 0$, i.e. $s_2 + 2s_3 = 0$. 
This gives $\gamma_k(\lam) = -s_4$, which we may assume by Case 1 is negative. Plugging in $k = 0$, we have
\begin{eqnarray*}
A &=&   - c m i^2\\
B &=& - d m i^2 - c n i^2 - 2 c m ij\\
C &=&     - d i^2n   - 2 d i jm  - 2 ci j  n - c j^2m\\
D &=&  - 2 d n i j - d m j^2 - c n j^2\\
E &=& -d n j^2.
\end{eqnarray*}
One can then check that $\disc_4(AU^4 + BU^3W + CU^2W^2 + DUW^3 + EW^4) = 0$. \qed

\subsection{Stable vectors}

\begin{Thm}\label{thm-4}
A vector $v \in \Vf$ is stable for the action of $\Gf$ if and only if $\Deltaf(v) \neq 0$. 
\end{Thm}

\textit{Proof}. As above, let $v = F_v(X, Y, Z) + (mU + nW)T$, with $F_v(X, Y, Z)$ defined as in (\ref{eqn-4v}), and define $A, B, C, D, E$ as in (\ref{eqn-4A-E}). By Proposition \ref{prop-hilmum} and Lemma \ref{lem-invar4}, if $\Deltaf(v) \neq 0$, then $v$ is stable.
Now suppose $\Deltaf(v) = 0$. Then $(mU + nW)\disc_{2}(F_v)$ has a double root, so after acting by $\SL_2 = \SL_{(U, V)}$ we may assume the coefficients of $U^4$ and $U^3W$ are zero,  i.e. $A = B = 0$.
Note that if we write $F_v$ as $U \otimes f_1(X, Y, Z) + W \otimes f_2(X, Y, Z)$, then $f_1 = aX^2 + cY^2 + eZ^2 + gXY + iXZ + kYZ$ and $A = m\disc_{2}(f_1)$. 

 If $\disc_2(f_1) \neq 0$, then $m = 0$ and $B = n(4ace + gki - eg^2 - ak^2 - ci^2) = n\disc_2(f_1)$, which implies $n = 0$. 
Then $v$ has no negative weights with respect to the cocharacter $\lam = -\check\omega_1$, so by Lemma \ref{lem-notstable}, $v$ is not stable.

Now suppose $\disc_{2}(f_1) = 0$. 
By \cite[Section 1.10]{Mumford}, after acting by $\SL_3 = \SL_{(X, Y, Z)}$ we can take $a = g = i = 0$ (note that because $\disc_{2, (X, Y, Z)}$ is invariant under the action of $\SL_3$, this doesn't change our assumption that $A = B = 0$). Then since $B = 0$ we have
\begin{eqnarray*}
B = bm(4ce - k^2) = 0.
\end{eqnarray*}
If $b = 0$, then $v$ has no negative weights with respect to the cocharacter given by taking $(s_1, s_2, s_3, s_4) = (-1, 2, -1, 0)$, with notation as in the proof of Lemma \ref{lem-invar4}. If $m = 0$, then $v$ has no negative weights with respect to the cocharacter given by $(s_1, s_2, s_3, s_4) = (-2, 2, -1, 0)$.
If $bm \neq 0$, then $eZ^2 + kYZ + cY^2$ has a double root. Now we may act by  
$\SL_{(Y, Z)} \subset \SL_3$ 
to obtain a vector with $c = k = 0$ (note that acting by $\SL_{(Y, Z)}$ does not change the fact that $a = g = i = 0$). This vector has no negative weights with respect to the cocharacter given by $(s_1, s_2, s_3, s_4) = (-2, 2, 0, -1)$.\qed

\begin{corollary}\label{cor-4}
\begin{enumerate}
\item[1.] The representation $\ZVf(\bF_p)$ contains stable vectors for all primes $p$.
\item[2.] Dual cyclic vectors are stable in $\ZVf(\bF_p)$ if and only if $p \neq 3$.
\end{enumerate}
\end{corollary}

\textit{Proof}.
A dual cyclic vector in $\ZVf(\bF_p)$ is conjugate to one of the form
$aUX^2 + dWY^2 + jWXZ + kUYZ + mUT$ for some $a, d, j, k, m \in \bF_p^\times$. 
We have $\Deltaf(aUX^2 + dWY^2 + jWXZ + kUYZ + mUT) = 
-27 a^2 d^4 k^4 j^8m^6$, which shows that dual cyclic vectors are stable if and only if $p \neq 3$. This gives part 2 of the corollary.
To see that stable vectors exist in $\ZVf(\FF_3)$, note that $\Deltaf(UX^2 + WY^2 + (U + W)Z^2 + WXY +  (U + W)YZ + UT) = 125$. \qed

\section{Stable 6-grading}\label{section-6}

For the 6-grading, defining an appropriate invariant polynomial is more delicate than in the other cases, so we first need to work with the representation of $\ZGs$ on the $\bZ$-module $\ZVs$, instead of working directly with the representation of $\Gs$ on $\Vs$. The normalized Kac coordinates of $y_6$ are $1 0 1 \Rightarrow 01$, so $\ZGs$ is generated by $\{\ZG_{\al_1}, \ZG_{\al_3}, \ZT\}$. By the same logic as in Section \ref{section-4}, we have $\ZG_{\al_1} \simeq \SL_2$, and since the $\al_3$-root string through $\al_4$ has length 2, we have 
 $\ZG_{\al_3} \simeq \SL_2$.
The group $\ZG_{\al_1} \times \ZG_{\al_3}$ is a subgroup of $\ZGs$, and as a representation of $\ZG_{\al_1} \times \ZG_{\al_3} \simeq \SL_2 \times \SL_2$, we have
\begin{equation*}
\ZVs \simeq (P_1(2) \otimes P_2(2)) \oplus (P_1(2) \otimes \triv) \oplus (\triv \otimes P_1(2)),
 \end{equation*}
where the first irreducible subrepresentation is generated by $f_{\al_2}$, the second is generated by $f_{-\al_0}$, and the third is generated by $f_{\al_4}$.

\subsection{An invariant polynomial}

Consider the representation $P_1(X, Y, Z) \otimes P_2(T, U, W)$ of $\SL_3 \times \SL_3$ (this is the representation corresponding to the stable 3-grading of $F_4$. See Section \ref{section-3}). We have that $\ZG_{\al_1} \times \ZG_{\al_3} \hookrightarrow \SL_3 \times \SL_3$ by mapping $\ZG_{\al_1}$ to $\SL_{(X, Y)}\subset \SL_{(X, Y, Z)}$, and mapping $\ZG_{\al_3}$ to $\SL_{(T, U)} \subset \SL_{(T, U, W)}$. Then 
$\ZVs$ is isomorphic to the $(\SL_{(X, Y)} \times \SL_{(T, U)})$-subrepresentation of $P_1(X, Y, Z) \otimes P_2(T, U, W)$ spanned by vectors of the form
\begin{equation}\label{eqn-vector6}
v = (aX + bY)\otimes T^2 + (dX + eY) \otimes U^2 + (gX + hY)\otimes W^2 + (jX + kY)\otimes TU + Z\otimes (oT + rU)W  
\end{equation}
for $a, b, d, e, g, h, j, k, o, r \in \ZZ$ (here for consistency we are labeling the coefficients as in Section \ref{section-3}). Thinking of $v$ as a degree-2 homogeneous polynomial in $T, U, W$, we may take the discriminat $\disc_{2, (T, U, W)}$, which gives a map $P_1(X, Y, Z) \otimes P_2(T, U, W) \to P_3(X, Y, Z)$ that is invariant under the action of $\SL_{(T, U)} \subset \SL_{(T, U, W)}$ and equivariant under the action of $\SL_{(X, Y)}$. Using the formula (\ref{eqn-disc-tern-quad}), we have
\begin{equation*}
\disc_{2, (T, U, W)}(v) = \be_1X^3 + \be_2X^2Y + \be_3XY^2 + \be_4Y^3 + \be_5XZ^2 + \be_6YZ^2
\end{equation*}
where
\begin{eqnarray*}
\be_1 &=& 4adg - gj^2\\
\be_2 &=& 4bdg + 4aeg + 4adh - hj^2 - 2gjk\\
\be_3 &=& 4beg + 4bdh + 4aeh - 2hjk - gk^2\\
\be_4 &=& 4beh - hk^2\\
\be_5 &=& -ar^2 - do^2 + ojr\\
\be_6 &=& -br^2 - eo^2 + okr.
\end{eqnarray*}
Let $F_1(X, Y) = \be_1X^3 + \be_2X^2Y + \be_3XY^2 + \be_4Y^3$, and let $F_2(X, Y) = \be_5X + \be_6Y$, so that
$\disc_{2, (T, U, W)}(v) = F_1(X, Y) + F_2(X, Y)Z^2$. Then $F_1F_2$ is a degree-$4$ homogeneous polynomial in the variables $X, Y$.
Explicitly, we have
\begin{equation*}
F_1F_2 = AX^4 + BX^3Y + CX^2Y^2 + DXY^3 + EY^4
\end{equation*}
where
\begin{eqnarray*}
A &=& \be_1\be_5\\
B &=& \be_1\be_6 + \be_2\be_5\\
C &=& \be_2\be_6 + \be_3\be_5\\
D &=& \be_3\be_6 + \be_4\be_5\\
E &=& \be_4\be_6.
\end{eqnarray*}
The discriminant $\disc_{4, (X, Y)}(F_1F_2)$ is a homogeneous polynomial over $\bZ$ in the variables $a, b, d, \dots, r$, and in fact $\disc_{4, (X, Y)}(F_1F_2) \in 16\bZ[a, b, d, \dots, r]$, so we set $\Deltas(v) = \frac{1}{16}\disc_{4, (X, Y)}(F_1F_2)$. Then $\Deltas$ is a degree-36 polynomial over $\ZZ$ in $a, b, d, e, \dots, r$ that is invariant under $\SL_{(X, Y)} \times \SL_{(T, U)}$. Because $\Deltas$ has coefficients in $\ZZ$, we may view it as a polynomial over $L$, and viewed in this way it is an invariant polynomial on $\Vs = \ZVs(F)$ for the action of $G_{\al_1} \times G_{\al_3}$. 

\begin{Rem}
One can view the embedding $\ZVs \hookrightarrow P_1(3) \otimes P_2(3)$ as coming from the natural inclusion $(\frakg_\bZ)_{x_6, 1} \hookrightarrow (\frakg_\bZ)_{x_3, 1}$, where $x_m = \frac{1}{m}\check\rho$ as in Section \ref{section-dual-cyclic}.
\end{Rem}

\begin{Prop}
The map $\Deltas$ is an invariant polynomial on $\Vs$ for the action of $\Gs$. 
\end{Prop}

\textit{Proof}. By the preceding discussion we have that $\Deltas$ is invariant under the actions of $G_{\al_1}$ and $G_{\al_3}$. 
Note that
$\Gs$ is generated by $\{G_{\al_1}, G_{\al_3}, \Tor_0\}$, where $\Tor_0$ is the subtorus of $\Tor$ with cocharacter group generated by $\{\check\omega_2, \check\omega_4\}$. Then $\Tor_0$ acts by the following weights on $\Vs$, where an $n$ in the row with irreducible subrepresentation $M$ and column $\check\omega_i$ indicates that $\check\omega_i(t)\cdot v = t^nv$ for all $v \in M$:
\begin{center}
\begin{tabular}{c|c|c}
& $\check\omega_2$ & $\check\omega_4$\\
\hline
$P_1 \otimes P_2$ & $-1$ & $0$\\
\hline
$P_1 \otimes 1$ & 3 & 2\\
\hline
$1 \otimes P_1$ & 0 & -1.
\end{tabular}
\end{center}
To prove the proposition, it suffices to show that $\Deltas$ is invariant under the actions of $\check\omega_2(t)$ and $\check\omega_4(t)$ for all $t \in L^\times$. By looking at the explicit formulas for the coefficients $\be_i$, we see that that $\check\omega_2(t)\cdot F_1 = tF_1$ and $\check\omega_2(t)\cdot F_2 = t^{-1}F_2$.
Similarly, $\check\omega_4(t)\cdot F_1 = t^2F_1$ and $\check\omega_4(t)\cdot F_2 = t^{-2}F_2$. The result follows.\qed

\begin{Lem}\label{lem-invar6}
The map $\Deltas$ satisfies the condition of Proposition \ref{prop-hilmum}.
\end{Lem}

\textit{Proof}. The proof is similar to those of Lemma \ref{lem-invar3} and Lemma \ref{lem-invar4}. As in those proofs, we let $\lam = s_1\check\omega_1 + s_2\check\omega_1 + s_3\check\omega_3 + s_4\check\omega_4$ be a nontrivial cocharacter. Then $\lam$ acts with the following weights with respect to the basis $\{X\otimes T^2, Y\otimes T^2, \dots, Z\otimes UW\}$ as in (\ref{eqn-vector6}): 
\begin{eqnarray*}
\gamma_a(\lam) &=& -s_2\\
\gamma_b(\lam)  &=& -s_1 -s_2\\
\gamma_d(\lam)  &=& -s_2 - 2s_3\\
\gamma_e(\lam)  &=& -s_1 - s_2 - 2s_3\\
\gamma_g(\lam)  &=& 2s_1 + 3s_2 + 4s_3 + 2s_4\\
\gamma_h(\lam)  &=& s_1 + 3s_2 + 4s_3 + 2s_4\\
\gamma_j(\lam) &=& -s_2 -s_3\\
\gamma_k(\lam)  &=& -s_1 - s_2 - s_3\\
\gamma_o(\lam)  &=& -s_4\\
\gamma_r(\lam)  &=& -s_3 - s_4.
\end{eqnarray*}
By this we mean that $\lam(t)\cdot X\otimes T^2 = t^{\gamma_a(\lam)}X\otimes T^2$ and similarly for the other basis vectors, labeled as in (\ref{eqn-vector6}).
Let $v \in \Vs$ be defined as in (\ref{eqn-vector6}), and suppose $v \in (\Vs)_\lam$. We must show $\Deltas(v) = 0$. The Weyl group of $\Gs$ is $\la w_{\al_1} \ra \times \la w_{\al_3} \ra \simeq \bZ/2\bZ \times \bZ/2\bZ$. Since $\al_i(\lam) = s_i = -(w_{\al_i}(\al_i))(\lam)$, using Lemma \ref{lem-weyl}, we can, and will, assume that $s_1 \geq 0$ and $s_3 \geq 0$ throughout the proof. 

Note that if $s_4 > 0$, then $\gamma_o(\lam) < 0$ and $\gamma_r(\lam) < 0$. Since $v \in (\Vs)_\lam$, this implies $o = r = 0$, and one can check that this implies $\Deltas(v) = 0$. So it suffices to assume that $s_4 \leq 0$. Similarly, if $s_2 > 0$, then $b = e = k = 0$, so 
$\Deltas(v) = 0$. So it suffices to assume $s_2 \leq 0$. Lastly,
if $\gamma_g(\lam) < 0$, then since $\gamma_g(\lam) \geq \gamma_h(\lam)$, we have $g = h = 0$, which implies $\Deltas(v) = 0$.
Thus it suffices to assume that $\gamma_g(\lam) \geq 0$. 
For easy reference, we list these assumptions we have just made:
\begin{itemize}
\item[B1.] $s_1 \geq 0$
\item[B2.] $s_2 \leq 0$
\item[B3.] $s_3 \geq 0$
\item[B4.] $s_4 \leq 0$
\item[B5.] $\gamma_g(\lam) \geq 0$.
\end{itemize}
We now consider two special cases.

Case 1: $\gamma_a(\lam) = 0$, i.e. $s_2 = 0$.

In this case, we have
\begin{eqnarray*}
\gamma_b(\lam) &=& -s_1\\
\gamma_d(\lam) &=&  - 2s_3\\
\gamma_e(\lam) &=& -s_1 - 2s_3\\
\gamma_g(\lam) &=& 2s_1 + 4s_3 + 2s_4\\
\gamma_h(\lam) &=& s_1  + 4s_3 + 2s_4\\
\gamma_j(\lam) &=& -s_2 -s_3\\
\gamma_k(\lam) &=& -s_1  - s_3.
\end{eqnarray*}
First suppose $s_3 = 0$. 
If $s_1 > 0$, then $\gamma_b(\lam) < 0, \gamma_e(\lam) < 0$, and $\gamma_k(\lam) < 0$, so $b = e = k = 0$, and $\Deltas(v) = 0$. If $s_1 = 0$, then 
since $\lam$ is nontrivial, we must have $s_4 \neq 0$, but this implies $\gamma_g(\lam) = 2s_4 < 0$, contradicting assumption B5.
Otherwise $s_3 > 0$, so $d = e = j = 0$ and $\Deltas(v) = 0$. This completes Case 1.

Case 2: $\gamma_b(\lam) \leq 0$.

If $\gamma_b(\lam) < 0$, then we must also have $\gamma_e(\lam) < 0$ and $\gamma_k(\lam) < 0$, so $b = e = k = 0$, so $\Deltas(v) = 0$. 
So we may assume $\gamma_b(\lam) = 0$. We have
\begin{eqnarray*}
\gamma_d &=& s_1 - 2s_3\\
\gamma_e &=&  - 2s_3\\
\gamma_g &=& -s_1 + 4s_3 + 2s_4\\
\gamma_h &=& -2s_1 + 4s_3 + 2s_4\\
\gamma_k &=& - s_3\\
\gamma_r &=& -s_3 - s_4.
\end{eqnarray*}
If $s_3 = 0$, then since $\lam$ is nontrivial, assumptions B1 and B4 imply that $\gamma_g(\lam) < 0$, contradicting assumption B5. So we have that $s_3 > 0$, and so $e = k = 0$. If $s_4 = 0$, then $r = 0$, and $\Deltas(v) = 0$.
Otherwise $s_4 < 0$, which implies $2\gamma_d(\lam) + \gamma_h(\lam) = 2s_4 < 0$, so at least one of $\gamma_d(\lam), \gamma_h(\lam)$ is negative. If $d = e = k = 0$, we have $\Deltas(v) = 0$, and if $e = k = h = 0$, we also have $\Deltas(v) = 0$, so $\Deltas(v) = 0$ in either case. This completes Case 2.

Now to finish the proof, note that by Case 2, we may assume that $\gamma_b(\lam) = -s_1 -s_2 > 0$, so 
$0 \leq \gamma_g < s_2 + 3s_3 + 2s_4 \leq 3(s_3 + s_4)$, so $s_3 + s_4 > 0$, so $\gamma_r(\lam) < 0$ and
$r = 0$. 
If $\gamma_d(\lam) < 0$, then $\gamma_e(\lam) < 0$, so $d = e = r = 0$ and $\Deltas(v) = 0$.
If $\gamma_d(\lam) \geq 0$, then $\gamma_h = -\gamma_b - 2\gamma_d + 2s_4 < 0$, so $h = 0$. By Case 1, we may assume $s_2 < 0$. If $\gamma_e(\lam) \geq 0$, then $\gamma_g(\lam) \leq s_2 + 2s_4 < 0$, contradicting our assumption B5. So $\gamma_e(\lam) < 0$ and $e = h = r = 0$. This implies $\Deltas(v) = 0$.\qed

\subsection{Stable vectors}

The next lemma will allow us to prove that if $\Deltas(v) = 0$, then $v$ is not stable.

\begin{Lem}\label{lem-notstable6}
Let 
\begin{equation*}
v = (aX + bY)\otimes T^2 + (dX + eY) \otimes U^2 + (gX + hY)\otimes W^2 + Z\otimes (oT + lU)W + (pX + qY)\otimes TU .
\end{equation*}
If any of the following hold, then $v$ is not stable:
\begin{enumerate}
\item[1.] $b = e = k = 0$
\item[2.] $d = e = k = 0$
\item[3.] $e = h = k = 0$
\item[4.] $e = k = r = 0$
\item[5.] $g = h = 0$
\item[6.] $o = r = 0$
\item[7.] $d = e = r = 0$
\item[8.] $e = h = r = 0$
\item[9.] $b = h = k = 0$
\item[10.] $d = j = r = 0$
\item[11.] $a = b = j = h = 0$.
\end{enumerate}
\end{Lem}

\textit{Proof}. 
By Lemma \ref{lem-notstable}, in each case it suffices to exhibit a cocharacter $\lam$ such that $\gamma_n(\lam) \geq 0$ for all $n \in \{a, b, d, e, g, h, j, k, l, r \}$. We do so in each case listed in the lemma. Here $(s_1, s_2, s_3, s_4)$ corresponds to the cocharacter $\lam = s_1\check\omega_1 + s_2\check\omega_2 + s_3\check\omega_3 + s_4\check\omega_4$.
\begin{enumerate}
\item[1.] $(1, 0, 0, 0)$
\item[2.] $(1, -1, 1, -1)$
\item[3.] $(2, -2, 1, -1)$
\item[4.] $(2, -2, 1, 0)$
\item[5.] $(-1, 0, 0, 0)$
\item[6.] $(0, 0, 0, 1)$
\item[7.] $(0, -1, 1, 0)$
\item[8.] $(1, -2, 1, 0)$
\item[9.] $(2, 0, -1, 0)$
\item[10.] $(-2, 0, 1, 0)$
\item[11.] $(-1, 2, -1, 0)$.
\end{enumerate}
\qed

\begin{Thm}\label{thm-6}
A vector $v \in \Vs$ is stable for the action of $\Gs$ if and only if $\Deltas(v) \neq 0$.
\end{Thm}

\textit{Proof}. 
Let $v \in \Vs$, write $v$ in the form given by (\ref{eqn-vector6}), and define $\be_1, \dots, \be_6$ as above.
By Proposition \ref{prop-hilmum} and Lemma \ref{lem-invar6}, we see that if $\Deltas(v) \neq 0$, then $v$ is stable, so it suffices to prove the converse. It will help to note that we may write $v$ as
\begin{equation}\label{eqn-v-var4}
v = X \otimes (aT^2 + jTU + dU^2) + Y \otimes (bT^2 + kTU + eU^2) + (gX + hY)\otimes W^2 + Z \otimes (oT + rU)W.
\end{equation}
We break into two cases, depending on the characteristic of $L$. 

Case 1: $\ch(L) \neq 2$.

Suppose $\Deltas(v) = 0$.  Then $F_1F_2$ has a double root. Since $\SL_{(X, Y)}$ acts transitively on lines in $X, Y$ (see, e.g., \cite[Lemma 2.6]{RomanoThesis}), we may assume without loss of generality that $X^2 \mid F_1F_2$. Since the coefficient of $Y^4$ in $F_1F_2$ is $\be_4\be_6$ and the coefficient of $XY^3$ is $\be_3\be_6 + \be_4\be_5$, we have that $\be_4\be_6 = \be_3\be_6 + \be_4\be_5= 0$. 

Note that 
$\disc_2 (bT^2 + kTU + eU^2) = k^2 - 4be$.
I first claim that if $k^2 - 4be = 0$, then $v$ is not stable. Indeed, in this case acting by $\SL_{(T, U)}$ we may assume $e = k= 0$, which implies $\be_4 = 0$, which in turn implies $\be_6\be_3 = 0$. Because $e = k = 0$, we have that $\be_6 = -br^2$ and $\be_3 = 4bdh$. 
Thus we are in one of the first four cases of Lemma \ref{lem-notstable6}, and $v$ is not stable.
Thus our claim is proven, and for the remainder of Case 1 we may assume $k^2 - 4be \neq 0$. 

Suppose $\be_6 \neq 0$. Then $\be_4 = 0$ and so $\be_3 = 0$. Note that
\begin{eqnarray*}
\be_3 &=& h(4bd + 4ae - 2jk) + g(4be - k^2) \text{ and }\\
\be_4 &=& h(4be - k^2).
\end{eqnarray*}
Since $k^2 - 4be \neq 0$, we have $g = h = 0$, and we are in case 5 of Lemma \ref{lem-notstable6}, so $v$ is not stable.

Thus we may assume $\be_6 = 0$, which implies $\be_4\be_5 = 0$. Acting by $\SL_{(T, U)}$, we may assume $r = 0$, so $\be_6 = -eo^2 = 0$ and $\be_5 = -do^2$. Since $\be_4\be_5 = 0$, we must have $dho = 0$. If $o = 0$, then we are in case 6 of Lemma \ref{lem-notstable6}. 
Otherwise $e = r = 0$ and $dh = 0$. If $d = 0$, then we are in case 7 of Lemma \ref{lem-notstable6}.
 If $h = 0$, then we are in case 8 of Lemma \ref{lem-notstable6}.
This completes Case 1.

Case 2: $\ch(L) = 2$

Suppose $\Deltas(v) = 0$. Acting by $\SL_{(X, Y)}$, we may assume that $h = 0$. Similarly, we may act by $\SL_{(T, U)}$ and assume $r = 0$. 
With these assumptions, we have that 
\begin{equation}\label{eqn-notstable-char2}
\Deltas(v) = e^2g^6o^{12}k^4(ej + dk)^4(b^2d^2 + a^2e^2 + bej^2 + bdjk + aejk + adk^2).
\end{equation}
If $o = 0$, we use Lemma \ref{lem-notstable6} case 6. If $e = 0$, we use Lemma \ref{lem-notstable6} case 8. If $g = 0$, we use Lemma \ref{lem-notstable6} case 5. Henceforth we assume $ego \neq 0$.
We can write $v$ in the form
\begin{equation*}
v = X \otimes(aT^2 + jTU + dU^2) + Y \otimes (bT^2 + kTU + eU^2) + gX \otimes W^2 + Z \otimes oTW.
\end{equation*}
Recall that for any $x \in L$, the element $z = \begin{psmallmatrix} 1 & x\\0 & 1\end{psmallmatrix} \in \SL_{(T, U)}(L)$ acts on $P_2(T, U)$ as $z \cdot G(T, U) = G(T, xT + U)$, so
\begin{equation*}
z\cdot v = X \otimes [(a + jx + dx^2)T^2 + jTU + dU^2] + Y \otimes [(b + kx + ex^2)T^2 + kTU + eU^2] + gX \otimes W^2 + Z \otimes oTW.
\end{equation*}
Similarly, the element $w = \begin{psmallmatrix} 1 & x\\0 & 1\end{psmallmatrix} \in \SL_{(X, Y)}(L)$ acts as
\begin{equation*}
w \cdot v = X \otimes [(a + bx)T^2 + (j + kx)TU + (d + ez)U^2] + Y \otimes (bT^2 + kTU + eU^2) + gX \otimes W^2 + Z \otimes oTW.
\end{equation*}

Now suppose $k = 0$, and let $z_0 = \begin{psmallmatrix} 1 & \sqrt{\frac{b}{e}}\\0 & 1\end{psmallmatrix} \in \SL_{(T, U)}(L)$. We have 
\begin{equation*}
z_0\cdot v = X \otimes((a + \frac{bd}{e})T^2 + jTU + dU^2) + Y \otimes eU^2 + gX \otimes W^2 + Z \otimes oTW,
\end{equation*}
so $z_0\cdot v$ is not stable by case 9 of Lemma \ref{lem-notstable6}. Henceforth we assume $k \neq 0$. 

If $ej + dk = 0$, then $ej = dk$, so $\frac{j}{k} = \frac{d}{e}$. Similar to the above, if we act by $w_0 = \begin{psmallmatrix} 1 & {\frac{j}{e}}\\0 & 1\end{psmallmatrix} \in \SL_{(X, Y)}(L)$, we have
\begin{equation*}
w_0 \cdot v = X \otimes(a + \frac{bj}{k})T^2  + Y \otimes (bT^2 + kTU + eU^2) + gX \otimes W^2 + Z \otimes oTW,
\end{equation*}
and $z_1 \cdot v$ is not stable by case 10 of Lemma \ref{lem-notstable6}.
Henceforth we assume $ej + dk \neq 0$, which implies by (\ref{eqn-notstable-char2}) that
\begin{equation}\label{eqn-char2b}
(bd + ae)^2 + (bj + ak)(ej + dk) = 0.
\end{equation}

Now because $k \neq 0$, we may still act by $w_0$ as defined in the previous paragraph, and we have
\begin{equation*}
w_0 \cdot v = X \otimes (a'T^2  + d'U^2) + Y \otimes (bT^2 + kTU + eU^2) + gX \otimes W^2 + Z \otimes oTW,
\end{equation*}
where $a' = a + \frac{bj}{k}$ and $d' = d + \frac{ej}{k}$. Since $ej + dk \neq 0$, we have $d' \neq 0$. Let $z_1 = \begin{psmallmatrix} 1 & \sqrt{\frac{a'}{d'}}\\ 0 & 1\end{psmallmatrix} \in \SL_{(T, U)}(L)$. Then 
\begin{equation*}
z_1w_0\cdot v =  X \otimes d'U^2 + Y \otimes (b'T^2 + kTU + eU^2) + gX \otimes W^2 + Z \otimes oTW,
\end{equation*}
where 
\begin{eqnarray*}
b' &=& b + k\sqrt{\frac{a'}{d'}} + \frac{a'e}{d'}\\
&=& b + k\sqrt{\frac{ak + bj}{dk + ej}} + \frac{e(ak + bj)}{dk + ej}.
\end{eqnarray*}
By (\ref{eqn-char2b}), we have that $bj + ak = \frac{(ae + bd)^2}{dk + ej}$, so we have
\begin{eqnarray*}
b' &=& b + \frac{k(ae + bd)}{dk + ej} + \frac{e(ak + bj)}{dk + ej},
\end{eqnarray*}
so
\begin{eqnarray*}
b'(dk + ej) &=& b(dk +ej) + k(ae + bd) + e(ak + bj)\\
&=& 0.
\end{eqnarray*}
Thus $b' = 0$, and $z_1w_0\cdot v$ is not stable by Lemma \ref{lem-notstable6} case 11.\qed

\begin{corollary}\label{cor-6}
\begin{enumerate}
\item[1.] For every prime integer $\mathsf{p}$, the representation $\ZVs(\bF_{\mathsf{p}})$ contains stable vectors.
\item[2.] Dual cyclic vectors in $\ZVs(\bF_{\mathsf{p}})$ are stable if and only if $\mathsf{p} \neq 2$.
\end{enumerate}
\end{corollary}

\textit{Proof}. A dual cyclic vector is conjugate to a vector of the form 
\begin{equation*}
v_c = bY \otimes T^2 + jX \otimes TU + eY \otimes U^2 + gX\otimes W^2 + rZ \otimes UW.
\end{equation*} 
Then $\Deltas(v_c) =16b^9 e^3 g^6 j^6 r^{12}$, so dual cyclic vectors are stable whenever $\ch(L) \neq 2$. To see that $\ZVs(\bF_2)$ contains stable vectors, note that if
\begin{equation*}
v = aX \otimes T^2 + dX  \otimes U^2 +  hY \otimes W^2 + (jX + kY)\otimes TU + Z\otimes  rUW, 
\end{equation*}
then $\Deltas(v) = 16 a^9 d^3 h^6  k^6r^{12} - 8 a^8 d^2 h^6 j^2 k^6 r^{12} + 
 a^7 d h^6j^4 k^6 r^{12} $. \qed

\section{Stable 8-grading}\label{section-8}

The normalized Kac Coordinates of $\frac{1}{8}\check\rho$ are $1~1~1 \Rightarrow 0~ 1$, so we see that $\Phi_8$ has type $A_1$, and $\Ge$ is generated by $\{\Tor, G_{\al_3}\}$. By the same logic as in the previous section, we see that $G_{\al_3} \simeq \SL_2$. By Lemma \ref{lem-SL2-reps}, as a representation of $\SL_2$, we have
\begin{equation*}
\Ve \simeq P_2(2) \oplus P_1(2) \oplus P_0(2) \oplus P_0(2),
\end{equation*}
where $P_2(2)$ is generated as a subrepresentation by $f_{\al_2}$, $P_1(2)$ is generated by $f_{\al_4}$, and the two trivial subrepresentations are spanned by $f_{\al_1}$ and $f_{-\al_0}$ respectively. We will write  $P_0(2)_{\al_1}$ and $P_0(2)_{-\al_0}$ when we need to distinguish between these two trivial subrepresentations.

A maximal torus of $G_{\al_3}$ has cocharacter group generated by $\check\al_3 = -2\check\omega_2 + 2\check\omega_3 - \check\omega_4$, so $G_{0, 8}$ is in fact generated by $\{\Tor_0, G_{\al_3}\}$, where $\Tor_0$ is the subtorus of $\Tor$ with cocharacter group generated by $\{\check\omega_1, \check\omega_2, \check\omega_4\}$. We have that $\Tor_0$ acts by the following weights, where a $j$ in row $P_k(n)$ and column $\check\omega_i$ indicates that $\check\omega_i(t)\cdot v = t^jv$ for all $v \in P_k(n)$:
\begin{center}
\begin{tabular}{ c | c | c | c }
 & $\check\omega_1$ & $\check\omega_2$ & $\check\omega_4$  \\
 \hline
 $P_2(2)$ & 0 & -1 & 0\\
 \hline
 $P_1(2)$ & 0 & 0 & -1\\
 \hline
 $P_0(2)_{\al_1}$ & -1 & 0 & 0\\
 \hline
 $P_0(2)_{-\al_0}$ & 2 & 3 & 2.
 \end{tabular}
\end{center}

 \subsection{An invariant polynomial}
 
For clarity of notation, let $Z = f_{\al_1}$ and $W = f_{-\al_0}$, so that each vector in $P_2(2) \oplus P_1(2) \oplus P_0(2) \oplus P_0(2)$ is of the form 
$(F_1(X, Y), F_2(X, Y), fZ, gW)$ for some $F_1 \in P_2(2) = P_2(X, Y), F_2 \in P_1(2) = P_1(X, Y)$, and $f, g \in L$. Consider the map $\Deltae$ given by 
 \begin{equation*}
 \Deltae((F_1, F_2, fZ, gW)) = f^4g^2\disc(F_1)\disc(F_1F_2).
 \end{equation*}
Explicitly, if $(F_1(X, Y), F_2(X, Y), fZ, gW) = (aX^2 + bXY + cY^2, dX + eY, fZ, gW)$, we have
\begin{equation*}
\Deltae((F_1, F_2, fZ, gW)) = (f^2g(b^2 - 4ac)(cd^2 - bde + ae^2))^2,
\end{equation*}
so $\Deltae$ is a degree-16 
homogeneous polynomial in $a, b, ..., g$.
\begin{Lem}
The map $\Deltae$ is invariant for the action of $\Ge$. 
\end{Lem}
\textit{Proof}. The map $\Deltae$ is clearly invariant for the action of $G_{\al_3} \simeq \SL_2$, since the map $P_2(2) \oplus P_1(2) \to P_3(2)$ given by $(F_1, F_2) \mapsto F_1F_2$ is $\SL_2$-equivariant and the discriminant polynomials are $\SL_2$-invariant. So we just need to check that $\Deltae$ is invariant for the action of $\Tor_0$. We have
\begin{eqnarray*}
\Deltae(\check\omega_1(t)\cdot (F_1, F_2, fZ, gW)) &=& \Deltae(F_1, F_2, tfZ, t^{-2}gZ)\\
&=& t^{-4}f^4t^{4}g^2\disc(F_1)\disc(F_1F_2)
\end{eqnarray*}
\begin{eqnarray*}
\Deltae(\check\omega_2(t)\cdot (F_1, F_2, fZ, gW)) &=& \Deltae(tF_1, F_2, fZ, t^{-3}gZ)\\
&=& f^4t^{6}g^2t^{-2}\disc(F_1)t^{-4}\disc(F_1F_2)
\end{eqnarray*}
\begin{eqnarray*}
\Deltae(\check\omega_4(t)\cdot (F_1, F_2, fZ, gW)) &=& \Deltae(F_1, tF_2, fZ, t^{-2}gZ)\\
&=& f^4t^{4}g^2\disc(F_1)t^{-4}\disc(F_1F_2),
\end{eqnarray*}
so $\Deltae$ is invariant for $\Ge$. \qed

\begin{Lem}\label{lem-invar8}
The invariant polynomial $\Deltae$ satisfies the hypothesis of Proposition \ref{prop-hilmum}.
\end{Lem}

\textit{Proof}. 
Let $\lam s_1\check\omega_1 + s_2\check\omega_2 + s_3\check\omega_3 + s_4\check\omega_4$ be a nontrivial cocharacter. 
Then $\lam$ acts
with the following weights 
with respect to the basis $\{X^2, XY, Y^2, X, Y, Z, W\}$, labeled by coefficients $a, b, c, d, e, f, g$ as above:
\begin{eqnarray*}
\gamma_a(\lam) &=& -s_2\\ 
\gamma_b(\lam) &=&  - s_2 -s_3\\ 
\gamma_c(\lam) &=&   -s_2 -2s_3\\ 
\gamma_d(\lam) &=& -s_4\\ 
\gamma_e(\lam) &=& -s_3 - s_4\\ 
\gamma_f (\lam) &=& -s_1\\
\gamma_g (\lam) &=& 2s_1 + 3s_2 + 4s_3 + 2s_4.
\end{eqnarray*}
Let $v = (aX^2 + bXY + cY^2, dX + eY, fZ, gW) \in (\Ve)_\lam$. We must show $\Deltae(v) = 0$. By Lemma \ref{lem-weyl}, we may replace $\lam$ with a Weyl-conjugate, so we may assume $s_3 \geq 0$.
If $\gamma_f(\lam) < 0$, by plugging $f = 0$ into the formula for $\Deltae$, we see that $\Deltae(v) = 0$. Similarly if $\gamma_g(\lam) < 0$, then $g = 0$ and we have $\Deltae(v) = 0$.
If $s_2 > 0$, then since $s_3 \geq 0$, we have $a = b = c = 0$, so $\Deltae(v) = 0$.
Lastly, if $s_4 > 0$, then since $s_3 \geq 0$, we have $d = e = 0$, so $\Deltae(v) = 0$. 
 Thus we may assume $s_1 \leq 0, s_2 \leq 0$
$s_4 \leq 0$, and $\gamma_g(\lam) \geq 0$.

Putting together these assumptions, we see that they imply $s_3 > 0$ (since $s_3 = 0$ implies $\gamma_g(\lam) < 0$). Note that if $s_2 = 0$, then the fact that $s_3 > 0$ implies $a = b = c = 0$, which again gives $\Deltae(v) = 0$, so we may further assume $s_2 < 0$. We then have $0 \leq \gamma_g(\lam) = 2s_1+ 3(s_2 + 2s_3) - 2s_3 + 2s_4 < 3(s_2 + 2s_3)$, so $\gamma_c(\lam) < 0$ and $c  = 0$. We also have $\gamma_g(\lam) = -2(\gamma_b(\lam) + \gamma_e(\lam) )+ s_2 < -2(\gamma_b(\lam) + \gamma_e(\lam))$, so either $b = 0$ or $e = 0$. In either case, we have $\Deltae(v) = 0$. \qed

\subsection{Stable vectors}

\begin{Thm}\label{thm-8}
A vector $v \in \Ve$ is stable for the action of $\Ge$ if and only if $\Deltae(v) \neq 0$.
\end{Thm}

\textit{Proof}. 
Let $v = (aX^2 + bXY + cY^2, dX + eY, fZ, gW) \in \Ve$.
By Proposition \ref{prop-hilmum} and Lemma \ref{lem-invar8}, if $\Deltae(v) \neq 0$, then $v$ is stable. 
Suppose $\Deltae(v) = 0$. If $g = 0$, then $v$ has no negative weights with respect to $-\check\omega_1$, so $v$ is not stable.
Similarly, if $f = 0$, then 
$v$ has no negative weights with respect to $\check\omega_1$, so $v$ is not stable.
So we may assume $fg \neq 0$, and thus either $b^2 - 4ac = 0$ or $cd^2 - bde + ae^2 = 0$. If the former, then $F_1$ has a double root, so acting by $\SL_2$ we can assume $a = b = 0$. Then $v$ has no negative weights with respect to cocharacter obtained by taking $(s_1, s_2, s_3, s_4) = (0, 2, -1, 0)$ with the notation of Lemma \ref{lem-invar8}, so $v$ is not stable. 

If $b^2 - 4ac \neq 0$, then $F_1$ is not a square, but $\disc(F_1F_2) = 0$, i.e. $F_1F_2$ has a double root. Since  $F_1$ is not a square, we must have 
$F_2 \mid F_1$. After acting by $\SL_2$, we may assume that $F_2 = Y$, so that $a = d = 0$. Then $(s_1, s_2, s_3, s_4) = (0, 1, -1, 1)$ yields a cocharacter with respect to which $v$ has no negative weights, and $v$ is not stable. \qed

\begin{corollary}\label{cor-8}
Stable vectors exist in $\ZVe(\FF_p)$ for all primes $p$.
\end{corollary}

\textit{Proof}. If $v = (XY, X + Y, Z, W)$, then $\Deltae(v) = 1$, so we see that in this case there are stable vectors in $\ZVe(\FF_p)$ for all $p$. \qed

\section{Applications to $p$-adic groups}\label{section-padic}

We now return to the notation of the introduction and prove the theorems stated there. First we sketch the construction of supercuspidal representations in \cite{ReederYu}, which motivated the results of this paper. See \cite[Section 2.5]{ReederYu} for details. As in the introduction, we let $\Gk$ be a split, absolutely simple connected reductive group over a nonarchimedean local field $k$ with residue field $\frakf$ of characteristic $p$, and $\bF$ is an algebraic closure of $\frakf$. We again write $\mathcal{A}$ for an apartment associated to a maximal split torus $\mathcal{T}$ of $\Gk$. As in the introduction, we choose a hyperspecial point $x_0 \in \mathcal{A}$, which allows us to identify $\mathcal{A}$ with $X_*(\mathcal{T}) \otimes \RR$, where $X_*(\mathcal{T})$ is the cocharacter group of $\mathcal{T}$. Fix a rational point $x \in \mathcal{A}$. For $r \geq 0$, we let $\Gk(k)_{x, r}$ be the Moy--Prasad filtration subgroup as defined in \cite[Section 3.2]{MP2} (in \cite{MP2} this subgroup is denoted $\mathscr{G}_{x, r}$).

Let $\Psi$ be the set of affine roots corresponding to $\mathcal{T}$, and let $r(x)$ be the smallest positive value in the set $\{ \psi(x) \mid \psi \in \Psi \}$. Then $\Gk(k)_{x, 0}/\Gk(k)_{x, r(x)}$  forms the ${\mathfrak{f}}$-points of a split reductive group $\Gk_x$ defined over $\mathfrak{f}$ (whose base change $(\Gk_x)_{\bF}$ is denoted $\mathsf{G}_x$ in \cite[Section 2.3]{ReederYu}). 
Given $r \geq 0$ we let $\Gk(k)_{x, r+} = \underset{s > r}\bigcup \Gk(k)_{x, s}$. 
We set 
\begin{equation*}
	\Vk_x : = \Gk(k)_{x, r(x)}/\Gk(k)_{x, r(x)+}.
	\end{equation*}
Then $\Vk_x$ is a finite-dimensional $\frakf$-vector space (denoted $\mathsf{V}_x(\frakf)$ in \cite[Section 2.5]{ReederYu}), and $\Gk_x$ has an action on $\Vk_x$ induced by conjugation in $\Gk(k)$. Recall that the representation of $\Gk_x$ on $\Vk_x$ can be seen as arising from a graded Lie algebra, where the grading is defined using $x$ as described in Section \ref{section-sl-reps}.

Given a stable vector  $\lam \in \check \Vk_x$ (as in the introduction, we call $\lam \in \check\Vk_x$ stable if it is stable after base change to $\bF$) and a nontrivial character $\chi: \mathfrak{f}^+ \to \CC^\times$, we consider the composition $\chi \circ \lam: \Vk_x \to \CC^\times$ as a character of $\Gk(k)_{x, r(x)}$ that is trivial on $\Gk(k)_{x,r(x)+}$. Then the compactly induced representation
	\begin{equation*}
	\pi_x(\lam) := \ind_{\Gk(k)_{x, r(x)}}^{\Gk(k)} (\chi \circ \lam) 
	\end{equation*}
	is a direct sum of finitely many irreducible supercuspidal representations of $\Gk(k)$ of depth $r(x)$ \cite[Propostion 2.4]{ReederYu}. These representations are said to be {epipelagic}, as defined in \cite{ReederYu}.

Note that the choice of $x$ and the orbit of $\lam$ give information about the resulting supercuspidal representations: the stabilizer of $\lam$ is used to determine the decomposition of $\pi_x(\lam)$ into irreducibles, as well as the formal degree of these irreducibles \cite[Proposition 2.4]{ReederYu}. The dimension of $\Gk_x$ is used to determine the Swan conductor of the Langlands parameters that should correspond to these representations under the local Langlands correspondence \cite[Section 7.1]{ReederYu}.

\begin{Rem}
For stable vectors $\lam, \lam' \in \check\Vk_x$, it is natural to expect that the irreducible subrepresentations of $\pi_x(\lam)$ and $\pi_x(\lam')$ are the same if and only if $\lam$ and $\lam'$ are in the same orbit, but this question was not explored in \cite{ReederYu} or in subsequent papers on the topic. I will return to this question in forthcoming work.
\end{Rem}

For the rest of the section we assume that $\Gk$ is of type $F_4$. We have the following theorem, as stated in the introduction.

\begin{Thm}
Suppose $\Gk$ is of type $F_4$ and $x$ is a rational point in $\mathcal{A}$. Assume that $x$ is not conjugate under the affine Weyl group to $x_0 + \frac{1}{2}\check\rho$.  
\begin{enumerate}
\item[1.]  If $\check{\Vk}_x(\bF)$ contains stable vectors, then so does $\check{\Vk}_x$.
\item[2.] There exists an invariant polynomial $\Delta_x$ on $\check{\Vk}_x(\bF)$ such that $v \in \check{\Vk}_x(\bF)$ is stable under the action of $(\Gk_x)_{\bF}$ if and only if $\Delta_x(v) \neq 0$.
\end{enumerate}
\end{Thm}

\textit{Proof}.
By \cite[Theorem 4.2]{FintzenRomano}, the representation of $(\Gk_x)_{\bF}$ on the base change $\check\Vk_x(\bF)$ contains stable vectors if and only if $x$ is conjugate under the affine Weyl group to $x_0 + \frac{1}{m}\check\rho$, where $m \in \{2, 3, 4, 6, 8, 12\}$. Thus we may restrict our attention to these $x$. The Coxeter number for $\Gk$ is $12$, and if $x$ is conjugate to $x_0 + \frac{1}{12}\check\rho$, then the statements in the theorem follow easily from \cite[Section 2.6]{ReederYu}. 

Let $H$ be a split, absolutely simple group  of type $F_4$ over $\frakf$, let $\mathfrak{h}$ be the Lie algebra of $H$, and write $Y_*$ for the cocharacter group of a maximal split torus of $H$. 
Using our choice of hyperspecial point $x_0 \in \mathcal{A}$, we may identify $Y_*$ with $\mathcal{A}$. 
Under this identification and with notation as in Section \ref{section-dual-cyclic}, given a rational point $x \in \mathcal{A}$, we have that $\Gk_x \simeq H_{x, 0}$ and under this isomorphism $\Vk_x \simeq \mathfrak{h}_{x, 1}$ (this follows from the proofs in \cite[Section 3.3.2]{RomanoThesis}). 
Thus describing stable vectors in $\check\Vk_x$ is equivalent to describing stable vectors in $\check{\mathfrak{h}}_{x, 1}$ under the action of $H_{x, 0}$. 
If $x$ is conjugate to $x_0 + \frac{1}{3}\check\rho$, then the claims in the theorem follow from Theorem \ref{thm-3} and Corollary \ref{cor-3}; if $x$ is conjugate to $x_0 + \frac{1}{4}\check\rho$, then the claims follow from Theorem \ref{thm-4} and Corollary \ref{cor-4}; if $x$ is conjugate to $x_0 + \frac{1}{6}\check\rho$, then the claims follow from Theorem \ref{thm-6} and Corollary \ref{cor-6}; and finally, if $x$ is conjugate to $x_0 + \frac{1}{8}\check\rho$, the claims follow from Theorem \ref{thm-8} and Corollary \ref{cor-8}. \qed

As a corollary, we have the following application to supercuspidal representations.

\begin{corollary}
Suppose $\Gk$ is of type $F_4$ and $x = x_0 + \frac{1}{m}\check\rho$ for some $m \in \{ 3, 4, 6, 8, 12\}$. Then $\check\Vk_x$ contains the input for the construction of \cite[Section 2.5]{ReederYu} as described above. 
Thus the construction of \cite{ReederYu} can be used to form epipelagic supercuspidal representations of $\Gk(k)$.
\end{corollary}

In particular, the corollary remains true if $k = \QQ_p$ and does not depend on $p$. We note that this is an improvement on \cite[Corollary 4.3]{FintzenRomano} since we obtain representations of $\Gk(k)$ without first moving to an unramified extension of $k$. 

Lastly, we revisit Conjecture \ref{conj}, which says, in the language of Section \ref{section-dual-cyclic}, that if $x$ is conjugate to $x_0 + \frac{1}{m}\check\rho$ and $m$ divides the Coxeter number of $\Gk$, then dual cyclic vectors are stable in $\check{\mathfrak{h}}_{x, 1}$ if and only if $\ch(\frakf) \nmid \frac{h}{m}$. 

\begin{Thm}
Conjecture \ref{conj} holds when $\Gk$ is of type $F_4$ for all $m \neq 2$.
\end{Thm}

\textit{Proof}. 
The conjecture follows from Corollary \ref{cor-3}, Corollary \ref{cor-4}, and Corollary \ref{cor-6}. 
\qed

\section{Appendix: Discriminant formulas}\label{section-disc}

In this appendix, we give explicit formulas for discriminant polynomials used throughout the paper. These formulas are classical (see, e.g., \cite{Salmon}, \cite{Sturmfels}), but in many classical texts, they are scaled in a way that makes them invalid in characteristics 2 and/or 3. In each case, we give references for the formulas used, which yield polynomials over $\bZ$ that are nonzero modulo $p$ for all primes $p$. 

\subsection{Binary quartics}

If 
\begin{equation*}
v = 
AX^4 + BX^3Y + CX^2Y^2 + DXY^3 + EY^4 \in P_4(2),
\end{equation*}
then
\begin{equation}\label{eqn-quarticdisc}
\begin{split}
\disc_{4}(v) =& ~B^2C^2D^2 - 4AC^3D^2 - 4B^3D^3 + 18ABCD^3 
-27 A^2 D^4 - 4 B^2 C^3 E + 16 A C^4 E\\& + 18 B^3 C D E - 
 80 ABC^2 DE - 6 AB^2 D^2 E + 144 A^2 C D^2 E - 27 B^4 E^2\\&  + 
 144 AB^2 CE^2 - 128 A^2 C^2 E^2 - 192 A^2 B D E^2 + 256 A^3 E^3.
 \end{split}
 \end{equation}
 This formula can be found explicitly in, for example, \cite[Section 4.1.1]{Cremona}.

\subsection{Ternary quadratics}

If
\begin{equation*}
v = AX^2 + BY^2 + CZ^2 + DXY + EXZ + FYZ \in P_2(3),
\end{equation*}
then
\begin{equation}\label{eqn-disc-tern-quad}
\disc_{2}(v) = 4ABC + DEF - AF^2 - BE^2 - CD^2.
\end{equation}
See, e.g., \cite[Section 1]{Lehman}.

\subsection{Ternary cubics}\label{section-disc-3-3}

The action of $\SL_3$ on $P_3(3)$ over an algebraically closed field $L$ has a degree-4 invariant generally denoted $S$, and a degree-6 invariant generally denoted $T$. Formulas may be found, e.g., in \cite[Sections 220, 221]{Salmon} or \cite[Sections 4.4, 4.5]{Sturmfels}), but we follow \cite{Fisher} for the scaling of the invariants. Explicitly,
let
\begin{equation*}
v = AX^3 + BX^2Y + CX^2Z + DXY^2 + EXYZ + FXZ^2 + GY^3 + HY^2Z + IYZ^2 + JZ^3,
\end{equation*}
and assume $A, B, \dots, J \in \bZ$. Set
\begin{eqnarray*}
S(v) &=&
E^4 - 8 DE^2 F + 16 D^2 F^2 + 24 CE F G - 48 B F^2 G - 8 CE^2 H -
 16 CD F H\\& & + 24 B EF H + 16 C^2 H^2 - 48 A F H^2 + 24 CDEI - 
 8 B E^2 I - 16 B D F I\\& & - 48 C^2 G I + 144 A F G I - 16 B C H I + 
 24 A E H I + 16 B^2 I^2 - 48 A DI^2\\& & - 48 CD^2 J + 24 B DEJ + 
 144 B C G J - 216 A E G J - 48 B^2 H J + 144 A DH J
\end{eqnarray*}
and
\begin{eqnarray*}
T(v) &=& -E^6 + 12 D E^4 F - 48 D^2 E^2 F^2 + 64 D^3 F^3 - 36 C E^3 F G + 
 144 C D E F^2 G\\& & + 72 B E^2 F^2 G - 288 B D F^3 G - 216 C^2 F^2 G^2 + 
 864 A F^3 G^2 + 12 C E^4 H\\ & &  - 24 C D E^2 F H - 36 B E^3 F H - 
 96 CD^2 F^2 H + 144 B D E F^2 H + 144 C^2 E F G H\\& & + 
 144 B CF^2 G H - 864 A E F^2 G H - 48 C^2 E^2 H^2 - 
 96 C^2 D F H^2 + 144 B C E F H^2\\& & + 72 A E^2 F H^2 - 
 216 B^2 F^2 H^2 + 576 A D F^2 H^2 + 64 C^3 H^3 - 288 A C F H^3\\& & - 
 36 C D E^3 i + 12 B E^4 I  + 144 C D^2 E F I - 24 B D E^2 F I - 
 96 B D^2 F^2 I\\& & + 72 C^2 E^2 G I + 144 C^2 D F G I - 
 720 B C E F G I + 648 A E^2 F G I + 576 B^2 F^2 G I\\& & - 
 864 A D F^2 G I + 144 C^2 D E H I - 24 B C E^2 H I - 36 A E^3 H I - 
 48 B C D F H I\\& & + 144 B^2 E F H I - 720 A D E F H I - 288 C^3 G H I + 
 1296 A C F G H I - 96 B C^2 H^2 I\\& & + 144 A C E H^2 I + 
 144 A B F H^2 I - 216 C^2 D^2 I^2 + 144 B C D E I^2 - 
 48 B^2 E^2 I^2\\& & + 72 A D E^2 I^2 - 96 B^2 D F I^2 + 576 A D^2 F I^2 + 
 576 B C^2 G I^2 - 864 A C E G I^2\\& & - 864 A B F G I^2 - 
 96 B^2 C H I^2 + 144 A C D H I^2 + 144 A B E H I^2 - 
 216 A^2 H^2 I^2\\& & + 64 B^3 I^3 - 288 A B D I^3 + 864 A^2 G I^3 + 
 72 C D^2 E^2 J - 36 B D E^3 J\\& & - 288 C D^3 F J + 144 B D^2 E F J - 
 864 C^2 D E G J + 648 B C E^2 G J - 540 A E^3 G J\\& & + 
 1296 B C d F G J - 864 B^2 E F G J + 1296 A D E F G J + 
 864 C^3 G^2 J\\& & - 3888 A C F G^2 J + 576 C^2 D^2 H J - 
 720 B CDE H J + 72 B^2 E^2 H J\\& & + 648 A DE^2 H J + 
 144 B^2 D F H J - 864 A D^2 F H J - 864 B C^2 G H J\\& & + 
 1296 A C E G H J + 1296 A B F G H J + 576 B^2 CH^2 J - 
 864 A C D H^2 J\\& & - 864 A B E H^2 J + 864 A^2 H^3 J + 
 144 B CD^2 I J + 144 B^2 D E I J - 864 A D^2 E I J\\& & - 
 864 B^2 C G I J + 1296 A C D G I J + 1296 A B E G I J - 
 288 B^3 H I J \\& &+ 1296 A B D H I J - 3888 A^2 G H I J - 
 216 B^2 D^2 J^2 + 864 A D^3 J^2 + 864 B^3 G J^2\\& & - 3888 A B D G J^2 + 
 5832 A^2 G^2 J^2.
\end{eqnarray*}
Then
\begin{equation}\label{eqn-disc-tern-cub}
\frac{1}{1728}(S^3 - T^2)
\end{equation}
is a homogeneous polynomial over $\bZ$ in the variables $A, B, \dots, J$, and the discriminant $\disc_3(v)$ is the image of this polynomial in $L[P_3(3)]$. See \cite[Section 4]{Fisher} for a more detailed discussion of this invariant, \cite[Section 7.2]{Fisher} for explicit formulas and scaling, and \cite[Section 10]{Fisher} for a more detailed look at the case when $\ch(L) \in \{2, 3\}$.

\bibliographystyle{plain}
\bibliography{f4stablevectors-arxiv}

\end{document}